# ON ITERATED LORENZ CURVES WITH APPLICATIONS*

ZVETAN IGNATOV[1], VILIMIR YORDANOV[2]

[1]*Sofia University “St. Kliment Ohridski”*
*Faculty of Economics and Business Administration*

[2]*Technical University Vienna*
*Financial and Actuarial Mathematics*

[2]*Vienna University of Economics and Business*
*Vienna Graduate School of Finance*

*e-mail:* vilimir.yordanov@tuwien.ac.at
vilimir.yordanov@vgsf.ac.at

*Zvetan Ignatov, Vilimir Yordanov.* ON ITERATED LORENZ CURVES WITH APPLICATIONS

It is well known that a Lorenz curve, derived from the distribution function of a random variable, can itself be viewed as a probability distribution function of a new random variable (*Arnold*, 2015). We prove the surprising result that a sequence of consecutive iterations of this map leads to a non-corner case convergence, independent of the initial random variable. In the primal case, both the limiting distribution and its parent follow a power-law distribution with exponent equal to the golden section. In the reflected case, the limiting distribution is the Kumaraswamy distribution with

* Zvetan Ignatov passed away on January 14, 2024. V.Y. would like to dedicate this paper to his memory. The paper is distributed following the ICMJE rules for authorship.

Zvetan Ignatov has worked on Section 1.1 and the conclusion. Vilimir Yordanov has worked on the introduction, Section 1.2, Section 2, the conclusion, and the appendix. Corresponding author. Non-academic e-mail: villyjord@gmail.com .

a conjugate value of the exponent, while the parent distribution is the classical Pareto distribution. Potential applications are also discussed.

## INTRODUCTION

The classical Lorenz curve finds numerous applications in applied statistics (*Arnold*, 2015), stochastic orders (*Arnold*, 2012), income inequality (*Kämpke* & *Radermacher*, 2015), risk analysis (Rüschendorf, 2013), portfolio theory (*Shalit* & *Yitzhaki*, 1984), etc. The pointed out popular references provide relevant details among many other existent excellent sources. The curve is a convenient tool for analysis of stochastic phenomena due to the plenty of important features it possesses as elaborated in (*Arnold*, 2015). However, a very distinguishable quality of it is the fact that the construct does not require a finite second moment of the parent distribution under focus. This constitutes an L1 norm approach to statistics, typical for the Italian school as discussed at length in (*Arnold*, 2015) and (*Hart*, 1975). Thus, it stays in contrast to the still currently predominant L2 norm approach, which directly relies on computing moments to a certain order in the process of characterization of the parent distribution. The effectively looser constraint which the Lorenz curve offers brings some benefits in general at the cost of higher technical difficulties encountered due to the need to work with differences of sample draws, and more generally, order statistics, instead of the centrality towards a moment. However, the L1 norm leads to more robustness, especially visible at present with the machine learning revolution, where algorithms based on such metrics are very popular as specified in (*Shalev-Shwartz* & *Ben-David*, 2014), (*Tibshirani*, 1996) and (*Hastie*, *Tibshirani*, *Friedman* & *Friedman*, 2009).

Within the above discussion, the Lorenz curve has also the very particular property that it can be viewed itself as a distribution function which provides certain modeling benefits. Instead of trying to target an ordinary characterization of the parent distribution just by heuristically inspecting the curve and constructing simple indices out of it, with the Gini being the most renown, we can focus entirely on studying the Lorenz curve itself. For example, moments can be computed from it as well as other relevant statistical quantities. They will both exist and by default behave in a more stable way. A thorough discussion about this approach can be found in (*Aaberge*, 2000) and (*Aaberge*, 2001). It is further reinforced by the fact that unlike the general distribution functions, the Lorenz curve has a [0, 1] domain and is also convex. This enables more precise numerical calculations and grants access to the powerful techniques of convex analysis.

With the above introduction, we not only highlight the benefits of using the Lorenz curve as a tool but also can distinguish a clear underlying iterative modeling logic. The natural question that occurs is how far we can go in it. Interestingly, such an inquiry somehow remained unexplored in the literature. We believe that it is worth the effort posing it since this can provide not only modeling insights but also would lead to interesting mathematical problems of standalone importance. For example, it is well known that a composition of Lorenz curves produces a Lorenz curve (e.g. (*Arnold*, 2015) and (*Sarabia*, 2008)), but what about a direct iteration. Such has already proven to be explored for other similar transforms with relevant implications. For example, in (*Harkness* & *Shantaram*, 1969) and (*Beek* & *Braat*, 1973), the stationary-excess distribution is under focus. Convergence results were proved with implications for reliability theory, point processes, Markov chains, insurance, etc. Even a generalization of classical mathematical results was made in (*Massey* & *Whitt*, 1993) and (*Lin*, 1994) in terms of coining a stochastic Taylor's theorem.

In this paper, we build on the logic of the last paragraph. Namely, we focus on Lorenz curve iterations. For the purpose, we formulate the problem explored in terms of several theorems and remarks. We provide the necessary proofs and discuss the implications that arise. Interestingly, to our knowledge, such natural line of research is completely missing in the literature and our paper sets the basis. The only exception is the study of (*Kämpke* & *Radermacher*, 2015), which seems to have been developed independently to our work and the preceding first preliminary public drafts of the latter[1]. In that book, iterated Lorenz curves are considered, however, the key convergence proof was not managed and was left as a challenging open problem (*Kämpke* & *Radermacher*, 2015, *Section* 7.2.1.1, p. 118). Here we are glad to provide the solution. Additionally, the iterated construct somehow remains artificial. No deep motivation is given why the aforementioned curves are considered and what their importance would be even within the income inequality setting which is the concrete theme of their research. In our opinion, they have big such with this being valid for different areas. With the current paper, we concentrate primarily on providing the technical toolkit, which has its own standalone mathematical implications. In the appendix, we take a small step toward practical applications by providing an illustrative example in quantitative portfolio management. The setting will be left for future research, focusing on technical generalizations and deeper applications.

[1] See the declaration at the end of the paper for the exact chronology.

## MAIN RESULTS

We formulate the results in terms of two theorems and several corollaries. Important remarks are also made.

### 1.1. PRIMAL CASE

**Theorem 1.** *Let X be an arbitrary non-negative random variable with a distribution function F and a positive finite mean $\mu_F$. It gives rise to a Lorenz curve, $L_F(x)$:*

$$L_F(x) = \frac{\int_0^x F^{-1}(u)\,du}{\mu_F} = \frac{\int_0^x F^{-1}(u)\,du}{\int_0^1 F^{-1}(u)\,du}, \tag{1}$$

*where $F^{-1}(u) = \inf\{y{:}F(y) \geq u\}$ for $0 \leq u \leq 1$ is the generalized inverse of $F(u)$ and $\mu_F = \int_0^1 F^{-1}(u)du < +\infty$ is the mean of X. We will denote by L both the Lorenz curve above (using the notations $L_F(x)$ and $L^F(x)$ interchangeably, emphasizing on the distribution function F that generates L and the same logic holds for $\mu_F$) and the operator $L(F(x))(.):[0,1] \to [0,1]$ which maps $F(.)$ to $L_F(.)$.*

*Since $L_F(x)$ by itself could be viewed as a distribution function of a random variable (with the possible extension $L_F(x) = 0$ for $x < 0$ and $L_F(x) = 1$ for $x > 0$ ), if for i = 0, 1, ... we consider the sequence of distribution functions $H_i^F(x)$:*

$$\begin{aligned}
&H_0^F(x) = F(x)\\
&H_1^F(x) = L^{H_0^F}(x) = L^F(x) = L_F(x) = L_0(x)\\
&H_2^F(x) = L^{H_1^F}(x) = L\left(L^F(x)\right) = L_1(x)\\
&H_3^F(x) = L^{H_2^F}(x) = L\left(L_1(x)\right) = L_2(x)\\
&\qquad\qquad \dots\\
&H_n^F(x) = L^{H_{n-1}^F}(x) = L\left(L_{n-2}(x)\right) = L_{n-1}(x)\\
&H_{n+1}^F(x) = L^{H_n^F}(x) = L\left(L_{n-1}(x)\right) = L_n(x),
\end{aligned} \tag{2}$$

*where both in H and L by the subscript n we indicate the iteration and by the superscript the starting distribution F, we have that $H_n^F(x)$ converges uniformly to the distribution function $G(x)$:*

$$G(x)=\begin{cases} x^{\frac{1+\sqrt{5}}{2}}, 0 \le x \le 1 \\ 0, x < 0 \\ 1, x > 1 \end{cases} \tag{3}$$

**Proof.** Different approaches can be employed for a proof. The most straightforward one is a resort to the theory of functional equations and using a variant of the contraction mapping theorem. However, we prefer a neater probabilistic solution. We provide a self-consistent proof in this direction.

From the definition of $L^F(x)$ in (1), we have the functional of the numerator to be continuous due to the integration and convex due to the monotonicity of the integrated cumulative distribution function (see for details – *Arnold*, 2015). Additionally, $L^F(x) = 0$ for $x = 0$ and $L^F(x) = 1$ for $x = 1$. Let's consider the functions:

$$\begin{aligned} &Z_0(x) = 0 \text{ for } 0 \le x \le 1 \\ &Z_0(1) = 1 \\ &Z_1(x) = x \text{ for } 0 \le x \le 1 \end{aligned} \tag{4}$$

It is easy to see that the following inequalities hold:

$$Z_0(x) \le L^F(x) \le Z_1(x) \text{ for } 0 \le x \le 1 \tag{5}$$

The left inequality is trivial, the right inequality follows from the convexity of $L_F(x)$. Taking the inverse functions, which preserve the monotonicity, and by suitable change of variables, we get:

$$Z_0^{-1}(x) \ge L_F^{-1}(x) \ge Z_1^{-1}(x) \text{ for } 0 \le x \le 1 \tag{6}$$

However, since taking inverses reverses the convexity in the particular monotonous case, we have that the functions $Z_0^{-1}(x), L_F^{-1}(x)$, and $Z_1^{-1}(x)$ are concave and additionally they pass through the points (0,0) and (1,1). So, we have also:

$$\begin{aligned} &Z_0^{-1}(x) = 1 \text{ for } 0 < x \le 1 \\ &Z_0^{-1}(0) = 0 \\ &Z_1^{-1}(x) = x \text{ for } 0 < x \le 1 \end{aligned} \tag{7}$$

Let now $x$ and $y$ be such that $0 \le x \le 1$ and $0 \le y \le 1$. If we consider the points (0,0) and ($y$,0), and take $x$ to be a scale parameter, the noted convexity of $L^F(x)$ gives us that the following is valid:

$$L_F(xy) \le xL_F(y) + (1-x)L_F(0) \tag{8}$$

$$L_F(xy) \le xL_F(y) \tag{9}$$

Analogously, the concavity of $L_F^{-1}(x)$ gives the validity of:

$$L_F^{-1}(xy) \ge xL_F^{-1}(y) + (1-x)L_F^{-1}(0) \tag{10}$$

$$L_F^{-1}(xy) \ge xL_F^{-1}(y) \tag{11}$$

Let's integrate both sides of (11) by *y* in the limits from 0 to *z*, where $0 \le z \le 1$. We have:

$$\int_0^z L_F^{-1}(xy)\,dy \ge \int_0^z xL_F^{-1}(y)\,dy \tag{12}$$

We can multiply the latter by *x* to get:

$$x\int_0^z L_F^{-1}(xy)\,dy \ge x^2\int_0^z L_F^{-1}(y)\,dy \tag{13}$$

Now, let's substitute $\upsilon = xy$ in the left integral of (13) and $\upsilon = y$ in the right one. We get:

$$\int_0^{xz} L_F^{-1}(v)\,dv \ge x^2\int_0^z L_F^{-1}(v)\,dv \tag{14}$$

Now, we can divide the two sides of the latter inequality by the positive constant $\int_0^1 L_F^{-1}(v)dv$ to get:

$$\frac{\int_0^{xz} L_F^{-1}(v)\,dv}{\int_0^1 L_F^{-1}(v)\,dv} \ge \frac{x^2\int_0^z L_F^{-1}(v)\,dv}{\int_0^1 L_F^{-1}(v)\,dv}, \tag{15}$$

which we can also write as:

$$L\left(L_F(v)\right)(xz) \ge x^2 L\left(L_F(v)\right)(z) \tag{16}$$

We can put $z = 1$ in (16) to get:

$$L\left(L_F\left(v\right)\right)\left(x\right) \geq x^2 \quad (17)$$

Now, take the inverse functions in the latter to get:

$$\left(L\left(L_F\left(v\right)\right)\left(x\right)\right)^{-1} \leq x^{\frac{1}{2}} \quad (18)$$

Let's introduce the following shorter notations:

$$L_0\left(x\right) = L_F\left(x\right), L_1\left(x\right) = L\left(L_F\left(x\right)\right), \ldots, L_n\left(x\right) = L\left(L_{n-1}\left(x\right)\right), \ldots \quad (19)$$

It is easy to see that the functions $L_n(x)$, $n = 0,1,2,...$, are convex and that:

$$L_n(0) = 0, n = 0,1,2,...,$$
$$L_n(1) = 1, n = 0,1,2,... \quad (20)$$

So, under the new notation, from (5), (9), (16), and (17), we have:

$$0 \leq L_0(x) \leq x \text{ for } 0 \leq x \leq 1 \quad (21)$$

$$L_0(xy) \leq xL_0(y) \text{ for } 0 \leq x \leq 1 \text{ and } 0 \leq y \leq 1 \quad (22)$$

$$x^2 \leq L_1(x) \text{ for } 0 \leq x \leq 1 \quad (23)$$

$$L_1(xy) \geq x^2 L_1(y) \text{ for } 0 \leq x \leq 1 \text{ and } 0 \leq y \leq 1 \quad (24)$$

If in (22) we replace $L_0(.)$ with $L_1(.)$, effectively jumping over one iteration, we can also get the inequality:

$$L_1(xy) \leq xL_1(y) \text{ for } 0 \leq x \leq 1 \text{ and } 0 \leq y \leq 1 \quad (25)$$

We will continue the proof by induction in terms of $n$.

For $n = 1$ we have from (24) and (25):

$$x^2 L_1(y) \leq L_1(xy) \leq xL_1(y) \text{ for } 0 \leq x \leq 1 \text{ and } 0 \leq y \leq 1 \quad (26)$$

Then we get consecutively:

$$x^2 L_2\left(y\right) \le L_2\left(xy\right) \le x^{\frac{3}{2}} L_2\left(y\right) \text{ for } 0 \le x \le 1 \text{ and } 0 \le y \le 1$$

$$x^{\frac{5}{3}} L_2\left(y\right) \le L_3\left(xy\right) \le x^{\frac{3}{2}} L_3\left(y\right) \text{ for } 0 \le x \le 1 \text{ and } 0 \le y \le 1 \tag{27}$$

and so on. Let's suppose now that all the inequalities are true till some $n$, i.e., let the last two inequalities for $n$ odd be:

$$x^{\alpha_{n+1}} L_n\left(y\right) \le L_n\left(xy\right) \le x^{\alpha_n} L_n\left(y\right) \text{ for } 0 \le x \le 1 \text{ and } 0 \le y \le 1 \tag{28}$$

and for $n$ even be:

$$x^{\alpha_n} L_n\left(y\right) \le L_n\left(xy\right) \le x^{\alpha_{n+1}} L_n\left(y\right) \text{ for } 0 \le x \le 1 \text{ and } 0 \le y \le 1 \tag{29}$$

where the sequence $\alpha_1 = 1$, $\alpha$, $= 2$, $\alpha_3 = \frac{3}{2}$, $\alpha_4 = \frac{5}{3}$, ... is given recursively by:

$$\alpha_{n+1} = 1 + \frac{1}{\alpha_n} \tag{30}$$

So, we suppose that for $n$ odd the inequalities (28) are valid, and we must prove that the following ones are valid:

$$x^{\alpha_{n+1}} L_{n+1}\left(y\right) \le L_{n+1}\left(xy\right) \le x^{\alpha_{n+2}} L_{n+1}\left(y\right) \text{ for } 0 \le x \le 1 \text{ and } 0 \le \text{y} \le 1 \tag{31}$$

Analogously, for $n$ even from the supposition that the inequalities (29) are true, we must prove that the following inequalities are true:

$$x^{\alpha_{n+2}} L_{n+1}\left(y\right) \le L_{n+1}\left(xy\right) \le x^{\alpha_{n+1}} L_{n+1}\left(y\right) \text{ for } 0 \le x \le 1 \text{ and } 0 \le \text{y} \le 1 \tag{32}$$

Let's start with the left inequality in (29). We have that the following holds:

$$x^{\alpha_n} L_n\left(y\right) \le L_n\left(xy\right) \text{ for } 0 \le x \le 1 \text{ and } 0 \le \text{y} \le 1 \tag{33}$$

Let's look at the two sides of the latter inequality as functions of $y$ with $x$ being just a parameter. From (33), for the inverse functions we have:

$$L_n^{-1}\left(\frac{y}{x^{\alpha_n}}\right) \ge \frac{\left(L_n^{-1}\left(y\right)\right)}{x} \text{ for } 0 \le x \le 1 \text{ and } 0 \le \text{y} \le 1 \tag{34}$$

In (34), we substitute $x = v^{\frac{1}{\alpha_n}}$ and $y = vw$ to get:

$$v^{\frac{1}{\alpha_n}} L_n^{-1}\left(w\right) \ge L_n^{-1}\left(vw\right) \text{ for at least } 0 \le v \le 1 \text{ and } 0 \le w \le 1 \tag{35}$$

Integrating the two sides of (35) by $w$ in limits from 0 to $z$ ($0 \le z \le 1$), we get:

$$v^{\frac{1}{\alpha_n}} \int_0^z L_n^{-1}(w)\,dw \ge \int_0^z L_n^{-1}(vw)\,dw \tag{36}$$

Multiplying both sides of (36) by $\upsilon$ and from the substitution $\tau = \upsilon w$ in the right integral, we get:

$$v^{1+\frac{1}{\alpha_n}} \int_0^z L_n^{-1}(w)\,dw \ge \int_0^{vz} L_n^{-1}(\tau)\,d\tau \tag{37}$$

Dividing the two sides of (37) by the positive constant $\int_0^1 L_n^{-1}(\tau)d\tau,$ we get:

$$\frac{v^{1+\frac{1}{\alpha_n}} \int_0^z L_n^{-1}(w)\,dw}{\int_0^1 L_n^{-1}(\tau)\,d\tau} \ge \frac{\int_0^{vz} L_n^{-1}(\tau)\,d\tau}{\int_0^1 L_n^{-1}(\tau)\,d\tau} \tag{38}$$

Due to (1) and (19), we can write the latter inequality as:

$$v^{\alpha_{n+1}} L_{n+1}(z) \ge L_{n+1}(vz) \tag{39}$$

for $0 \le \upsilon \le 1$ and $0 \le z \le 1$.

The inequality (39) coincides with the right inequality in (32). So, we proved that from the left inequality in (29) follows the right inequality in (32). Analogously, it follows that from the right inequality in (29) follows the left one in (32), and that from the left (right) inequality in (28) follows the right (left) inequality in (31). So, the induction is valid and for every integer the inequalities (31) and (32) are valid.

Now, let's substitute $y = 1$ in (31) and (32) and write, consecutively for $n = 1,2,...$, the inequalities that follow. We get for $0 \le x \le 1$:

$$0 \le F(x) \le 1$$

$$0 \le L_0(x) \le x$$

$$x^2 \le L_1(x) \le x$$

$$x^2 \le L_2(x) \le x^{\frac{3}{2}}$$

$$x^{\frac{5}{3}} \le L_3(x) \le x^{\frac{3}{2}}$$

$$x^{\frac{5}{3}} \le L_4(x) \le x^{\frac{7}{5}}$$

$$\begin{gathered} \ldots \\ x^{\alpha_{n+1}} \le L_n(x) \le x^{\alpha_n},\ n\ -\text{odd} \\ x^{\alpha_n} \le L_n(x) \le x^{\alpha_n+1},\ n\ -\text{even} \end{gathered} \tag{40}$$

Next, it is easy to see that the sequence (30) is convergent (bounded and monotonous in even and odd members, sharing a joint candidate limit). Solving for the joint positive limit, the equation:

$$\alpha = 1 + \frac{1}{\alpha} \tag{41}$$

gives us:

$$\lim_{n\to\infty} L_n(x) = x^{\frac{1+\sqrt{5}}{2}} \text{ for } 0 \le x \le 1 \tag{42}$$

The point-wise convergence of the sequence of functions $L_n(x)$ in the unit interval (compact set), together with their continuity and monotonicity in $x$, makes the continuity of the limiting function a necessary and sufficient condition for uniform convergence in the unit interval. Thus, the sequence of functions $L_n(x)$ are not only point-wise convergent but also uniformly convergent in [0,1]. This completes the proof for the primal case. ■

### 1.2. REFLECTED CASE

It is interesting to observe that the theorem can be adapted to a reflected case. For the purpose, we make first a remark and then formulate the result in terms of a consequent theorem.

**Remark 2.** The Lorenz curve, $L_F(x)$, can be reflected from the diagonal 1–x and thus it will give rise to a reflected Lorenz curve. Depending on how exactly the reflection is done, we can distinguish between two types of reflected Lorenz curves.

The first type can be viewed as 'simple-reflected'. It is defined by:

$$L_F^{s.ref}(x) = 1 - L_F^{-1}(1-x) \tag{43}$$

As prompted, the above expression means that we simply reflect the Lorenz curve from the diagonal. Direct differentiation gives that $L_F^{s.ref}(x)$ is increasing and convex with $L_F^{s.ref}(0) = 0$ and $L_F^{s.ref}(1) = 1$. This, together with the continu-ity of $L_F(x)$ and thus of $L_F^{s.ref}(x)$, gives that $L_F^{s.ref}(x)$ is indeed a Lorenz curve. The logic behind the architecture of (43) is straightforward. A problem with it is that we cannot define an iterated 'simple-reflected' Lorenz curve similarly to the way it was done in *Theorem 1* by a map. The only thing we can do is to apply consecutively the operator L from (2) starting from the initial distribution F and then make the reflection (43) of the resulted distribution. Surely, when the convergence in (2) holds, from the limit G(x), we will also get the limit:

$$G^{s.ref}(x)=1-G^{-1}(1-x)=\begin{cases}1-(1-x)^{\frac{\sqrt{5}-1}{2}},\ 0\le x\le 1\\ 0,\ x<0\\ 1,\ x>1\end{cases} \tag{44}$$

Thus, the iterative procedure lacks its own operator $L^{ref}(F(x))(.)$: $[0,1]\to[0,1]$ which to be applied in a standalone way at every step consecutively.

The second type has a bit more involved construct. It makes possible to form a reflected curve without an explicit reflection operation of a (primal) Lorenz curve. This allows to have the hinted above operator $L^{ref}(F(x))(.)$ which to be applied iteratively in a pure form. We will refer to this second curve as the 'reflected' curve, as it will serve as our baseline scenario. We reflect here rather the quantity $\frac{\int_{-\infty}^{x}(1-F(u))du}{\int_{-\infty}^{+\infty}(1-F(u))du}$, known as equilibrium distribution (integrated tail distribution) in insurance and reliability theory (*Klugman, Panjer* & *Willmot*, 2012), or under the alternative name of stationary-excess distribution, used mainly in renewal theory (*Whitt*, 1985). It[2] has plenty of applications in these areas, see (*Klugman, Panjer* & *Willmot*, 2012) and (*Embrechts, Klüppelberg* & *Mikosch*, 2013) for particular examples.

Concretely, we have the following theorem for the reflected Lorenz curve:

**Theorem 3.** *Let X be an arbitrary non-negative random variable with a distribution function F and a positive finite mean $\mu_F$. We can define a reflected Lorenz curve, $L_F^{ref}(x)$, by:*

$$L_F^{ref}(x)=1-\varphi_F^{-1}(x), \tag{45}$$

where:

[2] Note that it is not a Lorenz curve, so we do not have a reflection of a primal such curve. The reflection is indirect as indicated.

$$\varphi_F(x) = \frac{\int_{-\infty}^{1-x}(1-F(u))du}{\mu_F} = \frac{\int_{-\infty}^{1-x}(1-F(u))du}{\int_{-\infty}^{+\infty}(1-F(u))du} \tag{46}$$

*with* $\varphi_F^{-1}(u) = inf\{y\colon \varphi_F(y) \geq u\}$, $0 \leq u \leq 1$, *being the generalized inverse of* $\varphi_F(u)$ *and* $\mu_F = \int_{-\infty}^{+\infty}(1-F(u))du < +\infty$ *being the mean of X. We will denote by* $L^{ref}$ *both the Lorenz curve above (using the notations* $L_F^{ref}(x)$ *and* $L^{ref,F}(x)$ *interchangeably, emphasizing on the distribution function F that generates L) and the operator* $L^{ref}(F(x))(.)$: $[0,1] \to [0,1]$ *which maps* $F(.)$ *to* $L_F^{ref}(.)$. *The reflection is manifested by the upper integration limit in (46) and by the previous equation (45).*

*We claim that an adapted result to Theorem 1 holds also for the reflected Lorenz case:* $L_F^{ref}(x)$ *is a Lorenz curve, it can be viewed also as a distribution function (with the possible extension* $L_F^{ref}(x) = 0$ *for* $x < 0$ *and* $L_F^{ref}(x) = 1$ *for* $x > 0$*), and the corresponding limiting distribution after applying infinite* $L^{ref}(F(x))(.)$ *iterations on F is the Kumaraswamy distribution for which the primal parent is the classical Pareto distribution. Namely, if for* $i = 0,1,...$, *we consider the sequence of distribution functions* $H_i^{ref,F}(x)$, *we have:*

$$\begin{aligned}
&H_0^{ref,F}(x) = F(x) \\
&H_1^{ref,F}(x) = L^{H_0^{ref,F}}(x) = L^{ref,F}(x) = L_F^{ref}(x) = L_0^{ref}(x) \\
&H_2^{ref,F}(x) = L^{H_1^{ref,F}}(x) = L^{ref,F}(x) = L_F^{ref}(x) = L_1^{ref}(x) \\
&H_3^{ref,F}(x) = L^{ref,H_2^{ref,F}}(x) = L^{ref}\left(L_1^{ref}(x)\right) = L_2^{ref}(x) \\
&\qquad\qquad\qquad\qquad \ldots \\
&H_n^{ref,F}(x) = L^{ref,H_{n-1}^{ref,F}}(x) = L^{ref}\left(L_{n-2}^{ref}(x)\right) = L_{n-1}^{ref}(x) \\
&H_{n+1}^{ref,F}(x) = L^{ref,H_n^{ref,F}}(x) = L^{ref}\left(L_{n-1}^{ref}(x)\right) = L_n^{ref}(x)
\end{aligned} \tag{47}$$

*with* $H_n^{ref,F}(x)$ *converging uniformly to the distribution function* $G^{ref}(x)$, *where:*

$$G^{ref}(x) = \begin{cases} 1-(1-x)^{\frac{\sqrt{5}-1}{2}}, & 0 \leq x \leq 1 \\ 0, & x < 0 \\ 1, & x > 1 \end{cases} \tag{48}$$

**Proof.** The proof of the Lorenz case was given at length in *Theorem 1*. With suitable modification, the same logic could be applied to the reflected case. We will prove the result of the theorem by just elaborating on the main points.

First, by differentiation of $L_F^{ref}(x)$ at the definitions (45) and (46) and also by taking the values of $L_F^{ref}(x)$ at 0 and 1, it is clear that we have a Lorenz curve: $L_F^{ref}(x)$ is continuous, convex, with $L_F^{ref}(0)=0$, and $L_F^{ref}(1)=1$. The same holds for any function $L_n^{ref}(x)$ generated by the iterative procedure (47).

Second, we move to a counterpart of the induction procedure of *Theorem 1*. Suitable and essential modifications are needed in our case. The reasoning goes as follows.

Suppose that till some odd $n$ it is valid that:

$$1-(1-x)^{\frac{1}{\alpha_{n+1}}} \leq L_n^{ref}(x) \leq 1-(1-x)^{\frac{1}{\alpha_n}} \text{ for } 0 \leq x \leq 1, \tag{49}$$

or by noting that $\varphi_n(x) = L_n^{ref,-1}(1-x)$ (keeping to an obvious notation for the subscript iteration indices in line with the used one so far), and if we additionally introduce $\psi_n(x) = 1 - \varphi_n(x)$, in terms of $\varphi_n(x)$ and $\psi_n(x)$, we have:

$$1-(1-x)^{\alpha_n} \leq L_n^{ref,-1}(x) \leq 1-(1-x)^{\alpha_{n+1}} \tag{50}$$

$$1-x^{\alpha_n} \leq \varphi_n(x) \leq 1-x^{\alpha_{n+1}} \tag{51}$$

$$x^{\alpha_{n+1}} \leq \psi_n(x) \leq x^{\alpha_n} \tag{52}$$

Suppose also that till some even $n$ it is valid that:

$$1-(1-x)^{\frac{1}{\alpha_n}} \leq L_n^{ref}(x) \leq 1-(1-x)^{\frac{1}{\alpha_{n+1}}} \text{ for } 0 \leq x \leq 1 \tag{53}$$

or in terms of $\varphi_n(x)$ and $\psi_n(x)$:

$$1-x^{\alpha_{n+1}} \leq \varphi_n(x) \leq 1-x^{\alpha_n} \tag{54}$$

$$x^{\alpha_n} \leq \psi_n(x) \leq x^{\alpha_{n+1}}, \tag{55}$$

where the sequence $\alpha_1 = 1$, $\alpha_1 = 1$, $\alpha_3 = \frac{3}{2}$, $\alpha_4 = \frac{5}{3}$, ... is given recursively by:

$$\alpha_{n+1} = 1 + \frac{1}{\alpha_n} \tag{56}$$

We must prove that for odd $n$ it holds:

$$1-(1-x)^{\frac{1}{\alpha_{n+1}}} \leq L_{n+1}^{ref}(x) \leq 1-(1-x)^{\frac{1}{\alpha_{n+2}}} \text{ for } 0 \leq x \leq 1, \tag{57}$$

or in terms of $\varphi_{n+1}(x)$ and $\psi_{n+1}(x)$:

$$1-x^{\alpha_{n+2}} \le \varphi_{n+1}(x) \le 1-x^{\alpha_{n+1}} \tag{58}$$

$$x^{\alpha_{n+1}} \le \psi_{n+1}(x) \le x^{\alpha_{n+2}}, \tag{59}$$

and for even $n$ it holds:

$$1-(1-x)^{\frac{1}{\alpha_{n+2}}} \le L_{n+1}^{\text{ref}}(x) \le 1-(1-x)^{\frac{1}{\alpha_{n+1}}} \text{ for } 0 \le x \le 1, \tag{60}$$

or in terms of $\varphi_{n+1}(x)$ and $\psi_{n+1}(x)$:

$$1-x^{\alpha_{n+1}} \le \varphi_{n+1}(x) \le 1-x^{\alpha_{n+2}} \tag{61}$$

$$x^{\alpha_{n+2}} \le \psi_{n+1}(x) \le x^{\alpha_{n+1}} \tag{62}$$

It is more convenient to work with the $\varphi$*'s* due to the iteration formula (coming from the definition in (45) and (46) for $n \ge 0$):

$$\varphi_{n+1}(x) = \frac{\int_0^{1-x} \varphi_n^{-1}(u)\,du}{\int_0^1 \varphi_n^{-1}(u)\,du} \tag{63}$$

In terms of $\psi$*'s*, it translates to:

$$\psi_{n+1}(x) = \frac{\int_0^{x} \psi_n^{-1}(u)\,du}{\int_0^1 \psi_n^{-1}(u)\,du} \tag{64}$$

Now, we can apply the arguments from the proof of *Theorem 1*. Obviously, (64) has the form of (1). The implication (52) $\Rightarrow$ (59) is the same as the implication (28) $\Rightarrow$ (31), and the implication (55) $\Rightarrow$ (62) is the same as the implication (29) $\Rightarrow$ (32) by taking $y = 1$ in (28), (29), (31), and (32).

The statistical distribution from (48) is known in the literature as the (inflated) Kumaraswamy distribution (*Kumaraswamy*, 1980). It can be checked directly that its primal parent distribution (i.e., the resulting one after applying once to the obtained limit an inverse of the $L$ operator) is the classical[3] Pareto distribution with scale parameter equal to 1 and shape parameter equal to one plus the golden section (d.f.: $1-x^{-(1+\frac{1+\sqrt{5}}{2})}$). ■

[3] Or according to the standard classification presented in (*Arnold*, 2015, Section 4.2.8, p. 150) also Pareto (I) $\left(1, 1+\frac{1+\sqrt{5}}{2}\right)$.

**Remark 4.** It should be noted that although (*Kämpke* & *Radermacher*, 2015, Section 7.2.2.1, p. 119 and Section 7.2.2.2, p. 119 – 120) formulate a convergence hypothesis also for the reflected case, they again do not provide a proof. Additionally, it seems that they do not distinguish between the curves $L_F^{s.ref}(x)$ and $L_F^{ref}(x)$. It is not explicitly pointed out that the curves are not equal for arbitrary starting distribution. Also, no conceptual links to the useful equilibrium distribution were made.

## 2. DISCUSSION

Several important observations follow from the two theorems.

**Remark 5.** It should be mentioned that conceptually simpler but technically a bit messy proofs can be composed both for the primal and the reflected cases by just restricting the curves after the first iteration between 0 and x (here convexity is needed). Then taking inverses and integrating through the iterations without convexity arguments to be employed leads to a majorization by suitable polynomials. They are of the type $d_n\, x^{\gamma_n}$ with $d_n$ converging to 1 and $\gamma_n$ to the golden section for the Lorenz case and of the type $f_n + d_n\,(1 - x)^{\gamma_n}$ ) with $f_n$ and $d_n$ converging to 1 and $\gamma_n$ to the reciprocal of the golden section for the reflected Lorenz curve case.

**Remark 6.** It should be noted that both the primal operator, L(F(x))(.), and the reflected one, $L^{ref}(F(x))(.)$, act as two redistribution engines on the mass of the of the cumulative distribution function F(.). When applied iteratively to the generated d.f.'s from the Lorenz curves, the redistributions continue until convergence. Depending on the context, they can be viewed either as a set of income transfers, more formally Pigou-Dalton ones, e.g., (*Kämpke* & *Radermacher*, 2015) and (*Arnold*, 2015), or portfolio and risk rebalancings, e.g. (*Shalit* & *Yitzhaki*, 1984). More generally, such redistributions are well known in the theory of stochastic orders as discussed in (*Müller*, 2013) and (*Müller* & *Scarsini*, 2012). As seen from the corresponding formulas, the primal operator targets the left tail of the distribution, while the reflected one targets the right tail.

The latter observation can be seen alternatively also by considering the self-similarity differential equations from (*Kämpke* & *Radermacher*, 2015). Namely, we can first note that by using the upper and lower bound integral inequalities, for every Lorenz curve, it holds:

$$\frac{L(x)}{x} \le L'(x) \le \frac{1 - L(x)}{1 - x} \tag{65}$$

Now, relying on discrete intuition[4], we resort to the property of the Lorenz curves (see *Arnold* (2015) for details) that the atomic entities forming the distribution F (individual incomes, financial portfolio scenarios, risk scenarios, etc.) are proportional to the Lorenz density. So, (65) means that each atomic entity stays between the average of the lower and upper tail of the Lorenz curve.

If we impose L′ (x) to be a fraction $\varepsilon_u$ of the upper tail average (demonstrating upper self-similarity), we get the ODE:

$$L'(x) = \varepsilon_u \frac{1 - L(x)}{1 - x} \tag{66}$$

Its solution is:

$$L(x) = 1 - (1 - x)^{\varepsilon_u} \tag{67}$$

If we impose L′ (x) to be a multiple $\varepsilon_d$ of the down tail average (demonstrating down self-similarity), we get the ODE:

$$L'(x) = \varepsilon_d \frac{L(x)}{x} \tag{68}$$

Its solution is:

$$L(x) = x^{\varepsilon_d} \tag{69}$$

So, it becomes clear that the primal and the reflective iterative procedures produce self-similarity results as limit cases. Namely, for the former, each atomic entity becomes a fraction of the average of the upper tail of the Lorenz curve with a proportionality factor[5] $\varepsilon_u = \frac{\sqrt{5}-1}{2} \approx 0.618$. For the latter, each atomic entity becomes a multiple of the average of the down tail of the Lorenz curve with a proportionality factor $\varepsilon_d = \frac{\sqrt{5}+1}{2} \approx 1.618$.

**Remark 7.** It is interesting to observe how the two above-mentioned limits for the Lorenz curve and its reflected one can be obtained as solutions to fixed point problems.

---

[4] More precisely, we have that the d.f. *F* of *X* can be approximated by: $F_n = \frac{1}{n}\sum_{i=1}^{n} 1_{\{X_i \le x\}}$, where $X_i$ are i.i.d. draws (random sample) from *X* (Glivenko-Cantelli theorem). Observing that for the ordered sample $X_{i:n}$ it is valid $X_{i:n} = F^{-1}(U_{i:n})$ ($U_{i:n}$ – uniform ordered random sample), the Lorenz curve can be approximated by: $L\left(\frac{i}{n}\right) = \frac{\sum_{j=1}^{i} X_{j:n}}{\sum_{j=1}^{n} X_{j:n}}$. So, the mentioned proportionality is clear.

[5] Also named equity parameter in (*Kämpke* & *Radermacher*, 2015). .

Focusing only on Lorenz curve type of solutions[6], for the primal case, the limit G(x) solves the following functional equation:

$$\frac{\int_0^x G^{-1}(u)\,du}{\int_0^1 G^{-1}(u)\,du} = G(x) \tag{70}$$

After differentiation, it can be written alternatively as:

$$G'(x) = \frac{G^{-1}(x)}{\mu} \tag{71}$$

Now, we can make the substitution $g(x) = G^{-1}(x)$ and then apply the inverse function differentiation rule to it. After some manipulation, we get the following integro-differential equation:

$$g'(x) = \frac{\mu}{g(g(x))},\ \mu = g'(1) \tag{72}$$

with initial conditions (fixed points): $g(0) = 0$ and $g(1) = 1$. The book of (*Cheng & Li* (2008), Section 6.2.4, p. 212) provides details how to solve it in a closed form and get in the end the power-law for G(x). Actually, our iterated boundaries procedure gives as a by-product an alternative solution approach to the same equation, which in the spirit of the hinted possible general solution approach by the contraction mapping theorem should not be surprising.

For the reflected case, the limit $G^{ref}(x)$ solves the functional equation:

$$G^{ref}(x) = 1 - \varphi^{-1}(x)$$
$$\varphi(x) = \frac{\int_0^{1-x} (1 - G^{ref}(u))\,du}{\int_0^1 (1 - G^{ref}(u))\,du} \tag{73}$$

Again, with a suitable substitution (in this case simpler, namely directly $g_{ref}(x) = G^{ref}(x)$), it boils down to solving an integro-differential equation. Its form is:

$$g'_{ref}(x) = \frac{\mu}{g_{ref}(g_{ref}(x)) - 1},\ \mu = -g'(0), \tag{74}$$

with initial conditions (fixed points): $g_{ref}(0) = 0$ and $g_{ref}(1) = 1$. The book of

[6] As discussed in (*Kämpke*, 2011) there could be also other types having general properties, which will not be of interest to us due to our focus on Lorenz curves as a tool for stochastic analysis.

(*Cheng* & *Li* (2008), Section 6.2.4, Ex. 6.4, p. 214) again provides input on how to solve it in a closed form and get in the end the Kumaraswamy distribution. Our iterated boundaries procedure as a by-product again gives an alternative solution approach.

**Remark 8.** It is interesting to note that the functional equation (70) occasionally appears in the academic literature in one form or another, apart from the already pointed out instance in (*Cheng* & *Li*, 2008). For example, it can be traced back to (*McKiernan*, 1957)[7]. Additionally, as it was pointed out in (*Golomb, Marcus* & *Fine*, 1967), it has relation to the well-known Golomb sequence. Last but not least, it appeared twice in *The American Mathematical Monthly*: first as an ordinary problem in (*Nelson*, Aug. – Sep., 1968), with a solution provided in (*Hindmarsh*, Jun. – Jul., 1969), and then reappearing as a 'difficult problem' on a later occasion in (*Grossman*, Apr., 1976)[8].

Apart from the proofs and the remarks on the origin of the problem, it is interesting to see visually the convergence of the iterated primal, simple reflected, and reflected Lorenz curves. This is pointed out in Fig. 1, Fig. 2, and Fig. 3 respectively. In Fig. 4, the three iteration limits corresponding to the power and two times the Pareto laws are shown. The chosen starting distribution in all three cases is lognormal with mean of 0.5 and standard deviation of 0.2.

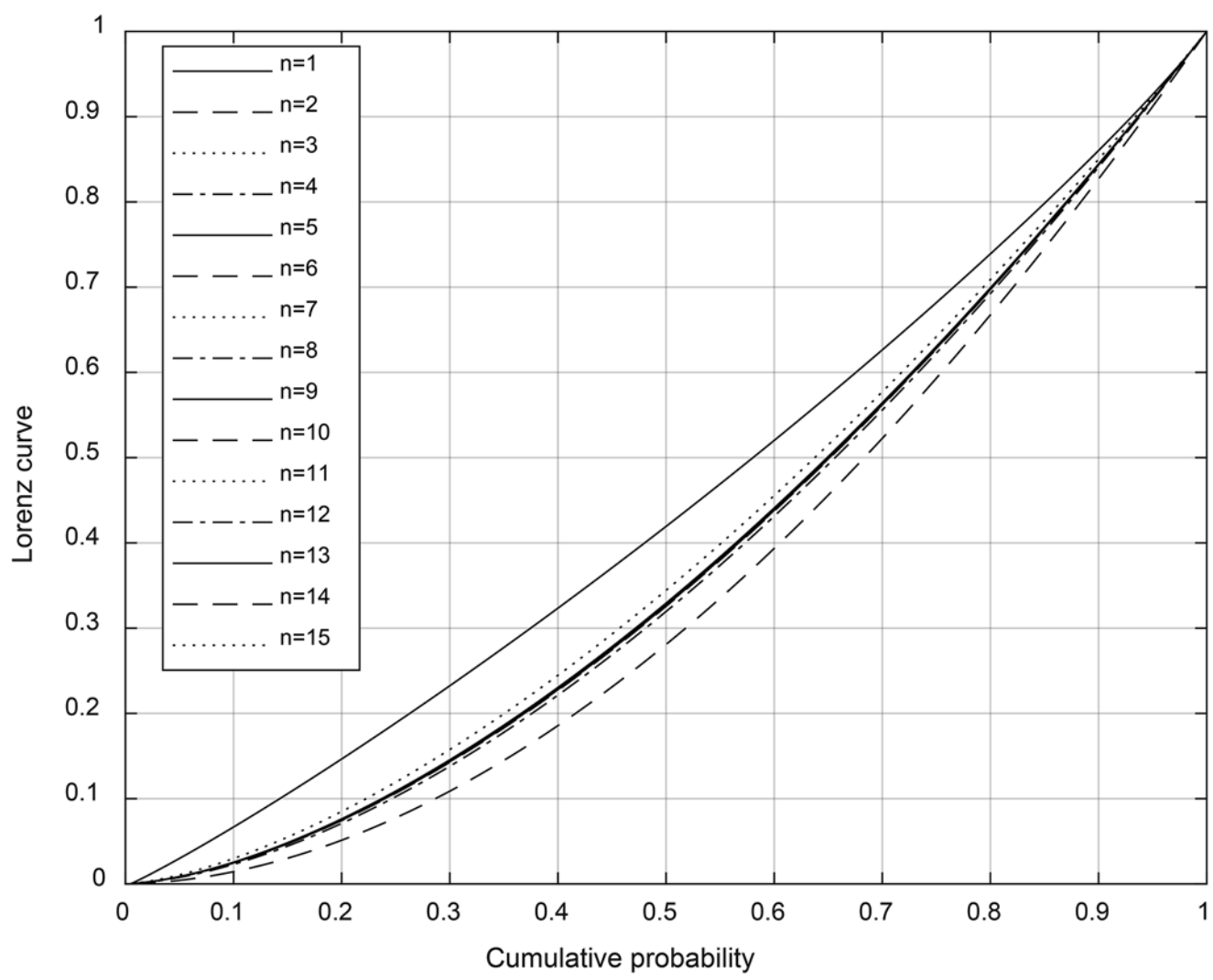

**Fig. 1**. Lorenz curve (primal) iterations

---

[7] In the form $f^{'}(x) = \dfrac{1}{f(f(x))}$.

[8] In the form $f'(x) = f^{-1}(x)$ with $f(0) = 0$ and $f'(x) > 0$ for $x > 0$.

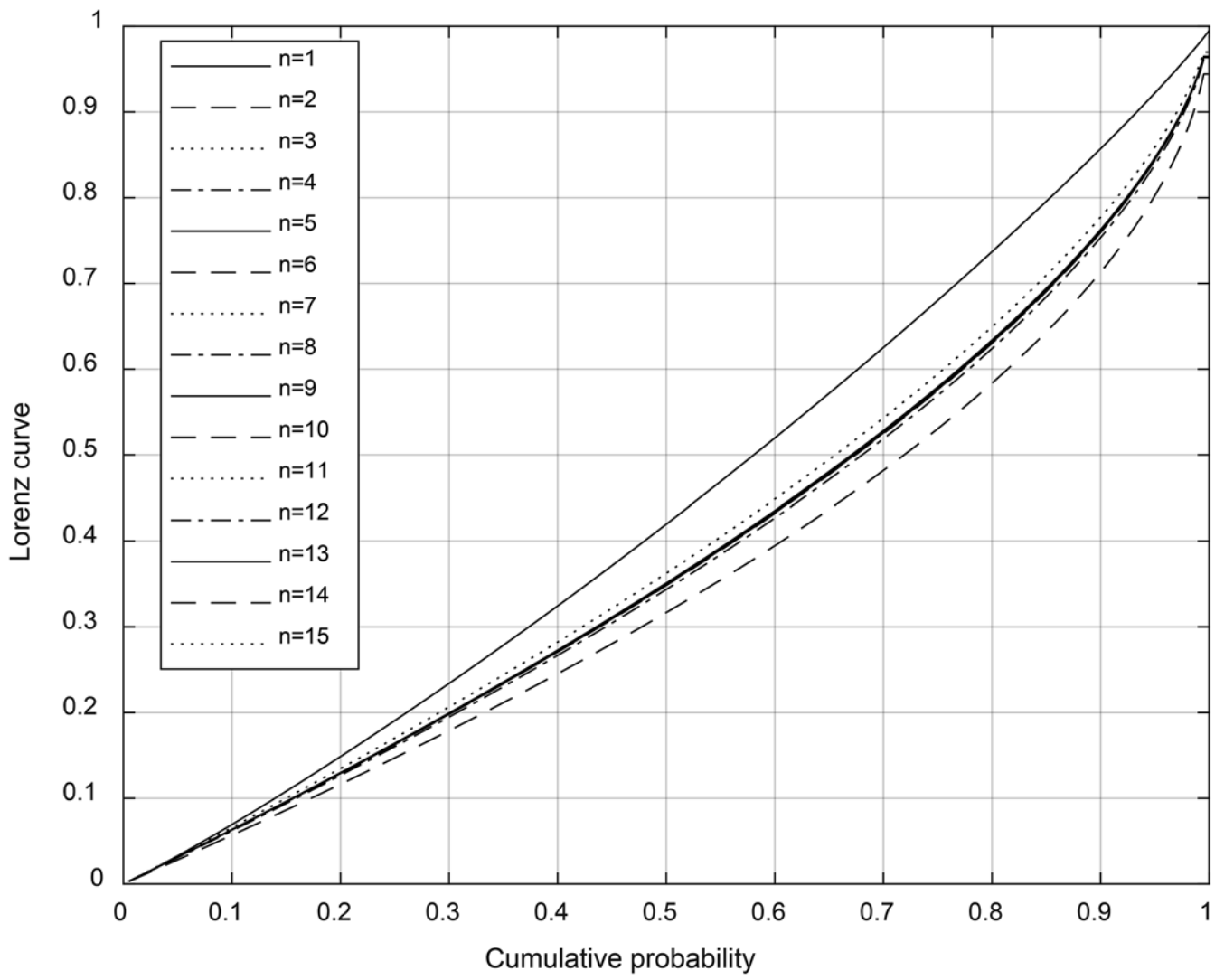

**Fig. 2**. Reflected Lorenz curve (simple reflection method) iterations

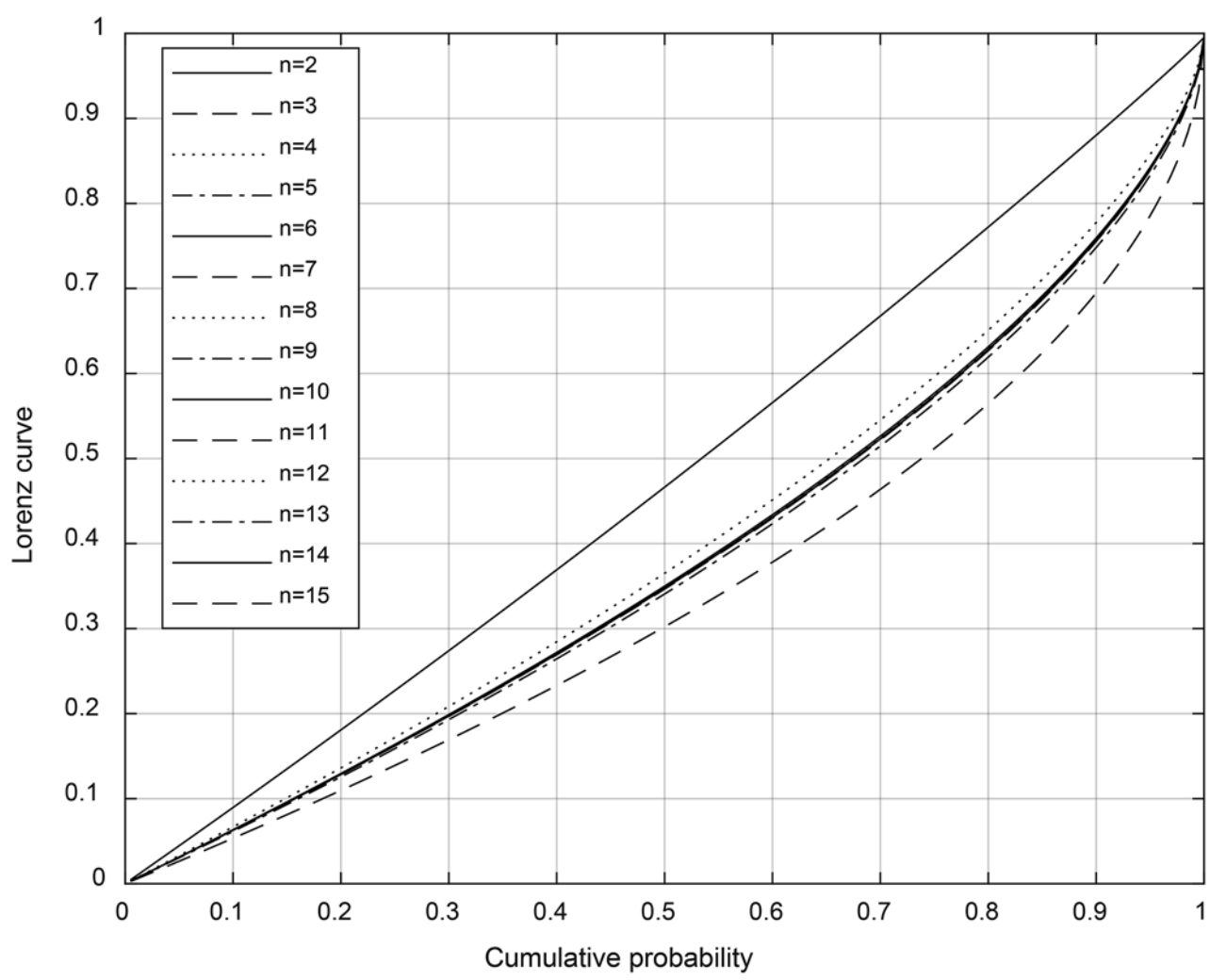

**Fig. 3.** Reflected Lorenz curve (iterated tail method) iterations

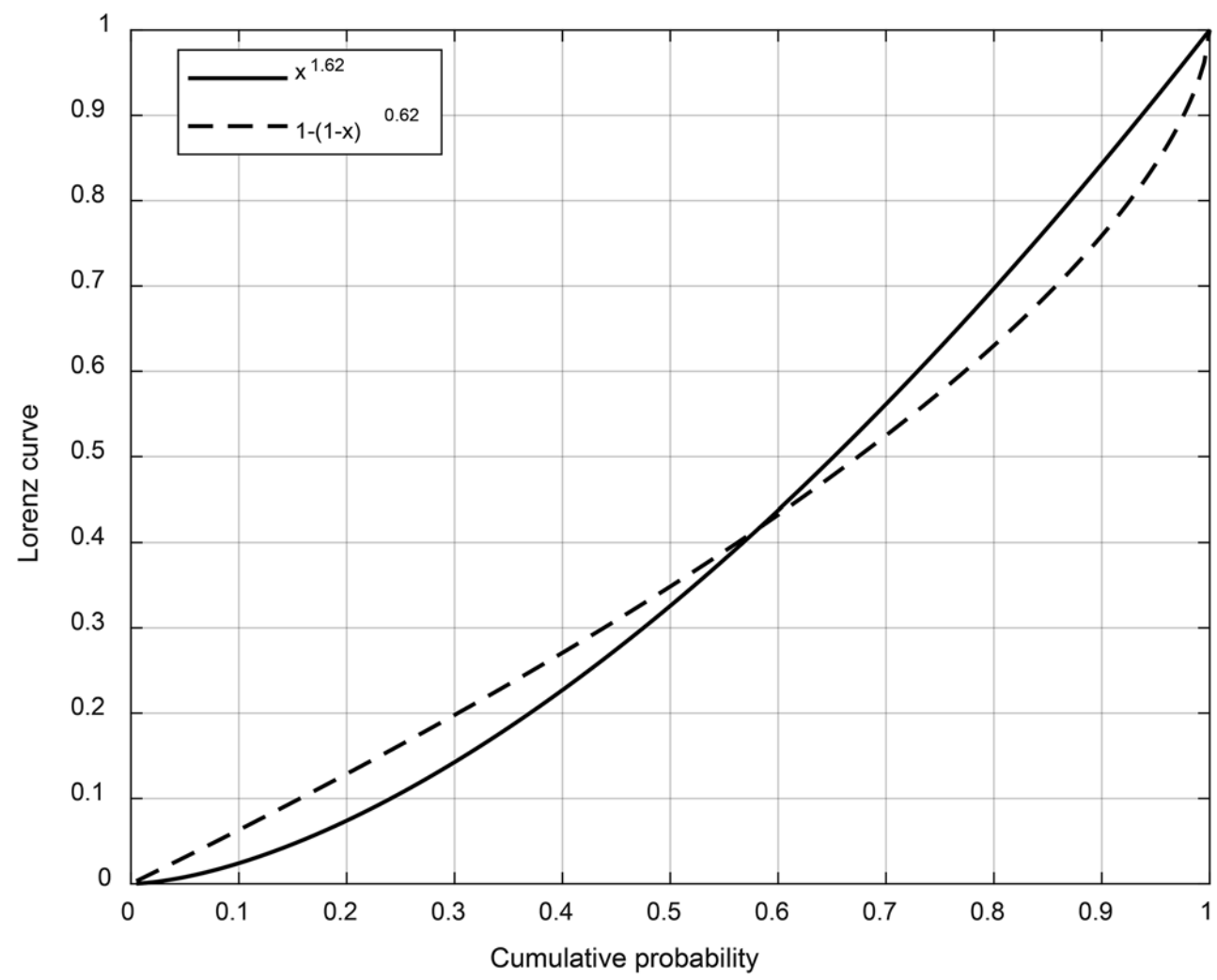

**Fig. 4**. Joint limits – primal, reflected (iterated tail method), and simple reflected Lorenz curves

We can observe very quick convergence in about 10 – 15 iterations. Our simulations[9] with different starting distributions do not seem to change the pattern.

## CONCLUSION

The presented results provide interesting and unexplored properties of Lorenz curves. Relevant scope for future work is open together with several possible applications in different directions and areas. From technical point of view, it would be interesting to link the convergence results to other similar known ones in the literature. E.g., we employed the equilibrium distribution and for it there are existent convergence properties as shown in (*Lin*, 1994). Of special interest would be the dual Lorenz curve $(L_F^{dual}(x) = 1 - L_F(1-x))$ to be explored under a similar iteration setting. If the primal curve gives the lower bound of a zonoid (see *Arnold* (2015) for details), the dual one gives the upper. A separate focused paper can also provide an alternative proof of the main results from the point of view of functional equations.

Apart from the pure mathematical relevance, the results of the paper can have multiple more practical applications. Such can be found in constructing new indices for dispersion in statistics, inequality in economics, and risk in quantitative finance. The two limit distributions also prompt connections to utility theory and decision under uncertainty.

---

[9] Source code: https://github.com/vjord-research/source-code/tree/main/lorenz-curves/paper1 .

**Declarations:**

**Chronology (ICMJE rules compatibility):** *Theorem 1* was stated for the first time in 2002 as a conjecture in the conference proceedings review article (*Chobanov* & *Ignatov*, 2002) on income inequality measures. No approach for a proof was shown. The current paper was started in 2004 – 2005 as a joint work when Vilimir Yordanov was assistant in Statistics to Zvetan Ignatov at Sofia University, Faculty of Economics and Business Administration targeting proofs, extensions, and applications. Vilimir Yordanov provided internally preliminary evidence for the validity of the conjecture in 2005 by numerical experiments and a technical proof in the spirit of *Remark 5*. Later that year, the more elegant proof from *Section 1.1* was accomplished by Zvetan Ignatov. Very preliminary drafts of the current paper were presented at: 1) Project conference "Les contraintes de l'economie Bulgare face à son adhésion à l'Union Européene" (funded by Rila project 2/9 – 2005 – 2007 as a cooperation between The Ministry of Education, Youth, and Science and The French Institute for financing grants and projects Egide) – proceedings in Bulgarian published in 2006 by the first author as an adapted version to the project's scope where he was a participant, 2) "The 2007 conference of the Society for Probability and Statistics of Romania". At both events, along with the side material, the second proof was announced. Then the paper had a consented status of an unpublished draft manuscript and was internally cited as such while getting updates by both authors through the years under irregular meetings' agenda. The final draft was agreed and composed in English in November 2022. The code realization, LaTeX typing, and typo corrections were completed in July 2023. The public release was agreed and scheduled for 2024. The current paper is a translation of the Latex version to MS Word and was announced on January 19, 2024. The final version was completed on November 12, 2024, after incorporating feedback from two referee reports and additional suggestions. These revisions led to minor technical corrections in the body text and the addition of an appendix.

## APPENDIX

### Introduction

The results developed in the body text, along with their standalone mathematical significance, can also have more applied implications. A good example would be the development of new statistical tools in economics and finance that rely on modified Lorenz curves and order statistics.

In (*Toshkov* & *Ignatov*, 2009), a paper citing the preliminary technical results of this work, it is observed that a weighted combination of the iteration limits of the primal and reflected Lorenz curves closely fits the Lorenz curve that characterizes income distribution in Bulgaria. The explanation given was that the economic dynamics of the country during the transition from a command to a market economy impacts quite unevenly the income inequality until it later stabilizes to a stationary level. As argued in Remark 6, the iterated Lorenz curves determine exactly such a 'stable' distribution pattern reached after probability mass redistributions. In economics of inequality context, they can be viewed as income transfers. Here we take another route. We will make an application to quantitative portfolio management. It serves as illustration of the tools developed and gives the major building blocks for further development[10].

In the context of portfolio theory, the probability mass redistributions can be interpreted as the result of rebalancing financial positions. Under a standard setting, an investor faces a universe of financial assets from which an optimal portfolio is formed through buying and selling. To achieve this goal, specific risk and return characteristics are targeted over a fixed time horizon. While taking the mathematical expectation is standard for estimating the latter parameter, the situation is different for the former. There is no consensus on how to measure risk. As discussed at length in (*Rachev*, *Stoyanov* & *Fabozzi*, 2008), different measures are possible. They largely depend on the specifics of the economic problem to be solved (e.g., the exact financial context within which the portfolio is formed – speculative, saving, pension, etc.; the type of the investor – institutional, bank, fund, insurance company, family office, sovereign fund, etc.; the overall balance sheet of the investor and the ALM strategy followed; legal & regulatory constraints; etc.) or the sole risk preferences of the investor, determined most generally by an abstract inherent utility function. The situation becomes even more complex, as each risk measure has its own specific mathematical qualities and characteristics[11]. Working with these measures requires a thorough understanding of them to assess the assets and limitations of the chosen approach. From this perspective, there are measures that seem to dominate others, though this is not absolute. As pointed out in (*Rachev*, *Ortobelli*, *Stoyanov*, *Fabozzi* & *Biglova*, 2008) and (*Rachev*, *Stoyanov* & *Fabozzi*, 2008), there are clearly desirable features to be sought. Apart from the classical coherence restriction introduced in (*Artzner*, *Delbaen*, *Eber* & *Heath*, 1999), the authors propose a set of others by providing almost an axiomatic treatment. The term an ideal risk measure is coined. In a flexible way, however, it is clearly elaborated that compromises can be made of some features over others depending on the context of the financial problem faced and the investor's attitude. We would like to position our analysis within such a setting. Namely, based on the exposition in the body text, we will propose a risk measure that does not belong to the class of ideal measures, but still addresses some important issues not encompassed by them. This paves the way for further and deeper theoretical analysis as well as more comprehensive empirical applications.

[10] We will postpone this to a later focused follow-up paper.

[11] Derived from the sound theory of probability distances, see for details (*Rachev*, *Klebanov*, *Stoyanov* & *Fabozzi*, 2013).

## 1. Modelling preliminaries

In general, risk measures try to provide a shortcut for characterizing the investor's attitude toward the statistical distribution of portfolio returns, so that an ordering of the investment outcomes can be achieved. In the classical Markowitz setting, the ranking is done solely based on the mean and the standard deviation (variance), with the latter being used as a risk measure. It is well known that the setting is inadequate. As discussed in (*Fabozzi* & *Pachamanova*, 2016), (*Rachev*, *Stoyanov* & *Fabozzi*, 2008), and (*Shalit* & *Yitzhaki*, 1984), among others, the problem is that characterizing the portfolio just by these two parameters is solely valid in the limiting cases of either normal (elliptical) distribution of the asset returns or quadratic investor's utility function, both of which do not hold in reality. So, investment decisions based on the mean-variance paradigm are biased and flawed.

A subtler issue is the general non-robustness of the optimal portfolios in this setting due to the ill-posedness of the optimization problem as investigated in (*Kan* & *Zhou*, 2007) and more recently also in (*Antonov*, *Lipton*, & *Lopez de Prado*, 2024). Such is mainly driven by parameters' instability in the inversion of large variance-covariance matrices. Even deeper, as argued in (*Rosenzweig*, 2021), the phenomenon of the winner-takes-all problem appears, meaning that the optimal portfolio is dominated by a small number of highly performing individual assets. Their impact prevents diversification which is certainly a mishap in contradiction to the major principles of financial theory. The partial fixes such as the ones to use the proposed by Markowitz himself semi-deviation measures (see *Markowitz*, 1959) or to resort to the risk parity approaches (see *Kolanovic* & *Wei*, 2013) do not support a fundamental way-out.

The most theoretically sound approach would be to get closer to the empirical facts about the asset returns together with the observed actual investors' behavior, which might be assumed to be driven by appropriate utility functions. The stochastic dominance rules that the investors obey to are supposed to be in line with the von Neumann–Morgenstern utility theory and its risk aversion implications. So, posing the validity of second order stochastic dominance (SSD) is a plausible assumption. If the Markowitz setting is characterized by mean-variance optimization, which is consistent with SSD only in specific cases, then under a more general framework, a new set of appropriate risk measures must be used to align fully with the prescribed stochastic dominance pattern. As argued in (*Levy* & *Robinson*, 2006), there are no absolute mean-risk measures that can guarantee the construction of universal SSD compliant efficient frontier in the optimization procedure. Every candidate measure would give only necessary conditions. Therefore, the choice should be carefully made so that it would take into consideration both the mathematical and the pure financial aspects of the problem. For example, the conditional value at risk measure (CVAR) elaborated in (*Rockafellar* & *Uryasev*, 2000) is not only a coherent measure but also is in line with SSD. Despite these positive features, its financial omnipotence is limited. The CVAR is primarily suitable for measurement of tail risks and is not appropriate for variability ones. The latter classification reflects the two universally accepted facets of financial risk as discussed in (*Righi* & *Ceretta*, 2015). The Gini's mean difference (GMD) introduced in (*Shalit* & *Yitzhaki*, 1984) as a risk measure also possesses coherence and SSD compliance. Since it can be viewed as a weighted average of CVARs (see *Mansini*, *Ogryczak* & *Speranza*, 2007), this makes it very suitable for measuring variability risks. It is important to note that, as already implied above, there are also additional qualities to be monitored, such as portfolio weights' stability and diversification effects[12]. The two aforementioned SSD consistent risk measures not always possess them. Considering all the above, it can be concluded that a risk measure should strike a reasonable balance between various conflicting features. They often enter a trade-off. We will target to construct a new risk measure that makes a reasonable such based on Lorenz curve applying the mathematical results from the main body.

[12] A much more general treatment on the desirable qualities of risk measures, especially in respect to classical risk theory, can be found in (*Cont, Deguest & Scandolo,* 2010) and (*Embrechts, Schied* & *Wang*, 2022). Considering them goes beyond the scope of the present paper.

## 2. Modelling details

We will motivate our risk measure construct on three pillars. First, as already mentioned, the developed iterated Lorenz curves methodology allows to view a specific sequence of portfolio rebalancing as an investment engine leading to a convergence to a stationary distribution. This operation allows to achieve a situation in which the investor won't be willing to further restructure. Thus, a 'stable' level of the portfolio's risk-return profile is attained[13].

Second, the natural question that arises is why this may be a desirable outcome to be pursued. Certainly, it should be related to the utility function of the investor and his preferences. Let's take a more careful look on how risk measures can be constructed from a Lorenz curve. The latter just orders the return scenarios as it can be seen from the numerator of the representation $L_F\left(\frac{i}{n}\right) = \frac{\sum_{j=1}^{i} X_{j:n}}{\sum_{j=1}^{n} X_{j:n}} = \frac{\sum_{j=1}^{i} X_{j:n}}{\mu_F}$ discussed before. Exactly this makes it a convenient tool for building risk measures. Weights can be put on the scenarios according to the risk aversion of the investor. This is exactly the case for the CVAR and the GMD. Namely, following (*Mansini*, *Ogryczak* & *Speranza*, 2007) and (*Mansini*, *Ogryczak* & *Speranza*, 2015)[14], after a discretization[15], the former can be written as:

$$CVAR_F(\alpha) = \mu_F \frac{L_F(1-\alpha)}{1-\alpha} \approx \sum_{j=1}^{n} \mathrm{w}_j^{CVAR} X_{j:n}, \tag{75}$$

where $\mathrm{w}_j^{CVAR} = -\frac{1}{\alpha n}$ for $j = 1, ..., i-1$ and $\mathrm{w}_i^{CVAR} = -1 + \sum_{j=1}^{i-1} \frac{1}{\alpha n}$ with $i = [an]$[16], and the latter as:

$$GMD_F = 4\mu_F \int_0^1 \left(\mathrm{x} - L_F(\mathrm{x})\right) \mathrm{dx} \approx \sum_{j=1}^{n} \mathrm{w}_j^{GMD} X_{j:n}, \tag{76}$$

where $\mathrm{w}_j^{GMD} = 2\left(\frac{2j-1-n}{n(n-1)}\right)$.

---

[13] The same logic is given in (*Toshkov* & *Ignatov*, 2009) within income inequality context.

[14] It must be noted that we use the market definition of CVAR which slightly differs from the one in the references.

[15] We won't consider here the Taleb's type of critique from (*Fontanari*, *Taleb* & *Cirillo*, 2018) about the need for a correction in the above discretization of the GMD. We will keep to the established practice. The biases that may arise will be considered in additional research.

[16] Here by [] we denote the ceiling function.

A large class of tail and variability measures as the ones above are dependent on a suitable Lorenz curve which acts as a base[17]. When the latter stabilizes by applying the iteration algorithms from *Theorem 1* or *Theorem 3*, we can expect to have reached a situation of minimum variability in the tails of the resulting distribution. This can be seen more directly by remembering the dual nature of a Lorenz curve. Namely, on one hand, it is an object characterizing the parent distribution, on the other hand, it is a distribution function on its own. So, in its latter apparition, the curve can itself be characterized by suitable risk measures. When the latter stop to differ gradually from the same ones applied to the parent distribution due to the Lorenz curves convergence, we may say that in the end of the iterative procedure, we have reached a situation of minimum variability in the tails. It is a completely plausible requirement the investor to have such a utility function that he is not willing to take tail risks of all orders. In case of standard statistical high moments existence, this investor would be averse to them. Compared to SSD, this is a more general requirement. Similar situation of high moments risk aversion is considered in (*Rosenzweig*, 2021) and has proved to be consistent with better diversification outcomes.

Third, the limiting stationary distributions from the iterative procedure seem also to be in line with the reality in financial markets. There are many studies in empirical asset pricing, e.g. (*Gabaix*, *Gopikrishnan, Plerou* & *Stanley*, 2003), showing that the returns of plenty of financial indices over short to medium horizons, as well as individual assets, follow power laws if very long data history is considered for estimation. So, the theoretically stationary distributions we derived seem also to empirically behave as such.

The above considerations motivate us to devise a risk measure in the spirit of GMD, where the investor, when constructing an optimal portfolio, aims for a distribution that mimics the tails of the limiting distributions while remaining close to SSD in the belly. This approach allows to have simultaneously minimum tails variability and SSD consistence both of which would reflect quite a reasonable investment behavior. The Lorenz curves of two possible targeted optimal parent distributions[18] are shown in Fig. 5. The first considers both the upper and lower

[17] A much more general treatment would be to consider the construction of the weighted ordered scenarios though Choquet integrals, with some basic building blocks of the applications of this method in mathematical risk theory given in (*Denuit*, *Dhaene, Goovaerts* & *Kaas*, 2006).

[18] It must be highlighted that a Lorenz curve can take negative values simply by construct if the parent random variable takes such. In this case, we speak about a generalized curve. It deviates from the classical definition used in the body text, which is limited to non-negative parent distributions. Under positive mean assumption (completely economically justified for optimal financial portfolios), we will have first a decreasing curve pattern then followed by an increasing one. Its intersection point with the x-axis corresponds to the value of the parent distribution's d.f. at naught. If the monotonicity of the classical case is broken, the convexity is theoretically preserved. Such Lorenz curve can be viewed as an extension of the classical definition and is well known in the liter-

tails, distinguishing between them. This is suitable for the general case where both positive and negative returns can occur. The second focuses on the absolute values of returns, thus not distinguishing between the two tails, with an emphasis on overall tail behavior. Both approaches depict also the non-stochastic belly as a straight line[19]. In each case, the modeler is allowed to have some flexibility up to scaling constants.

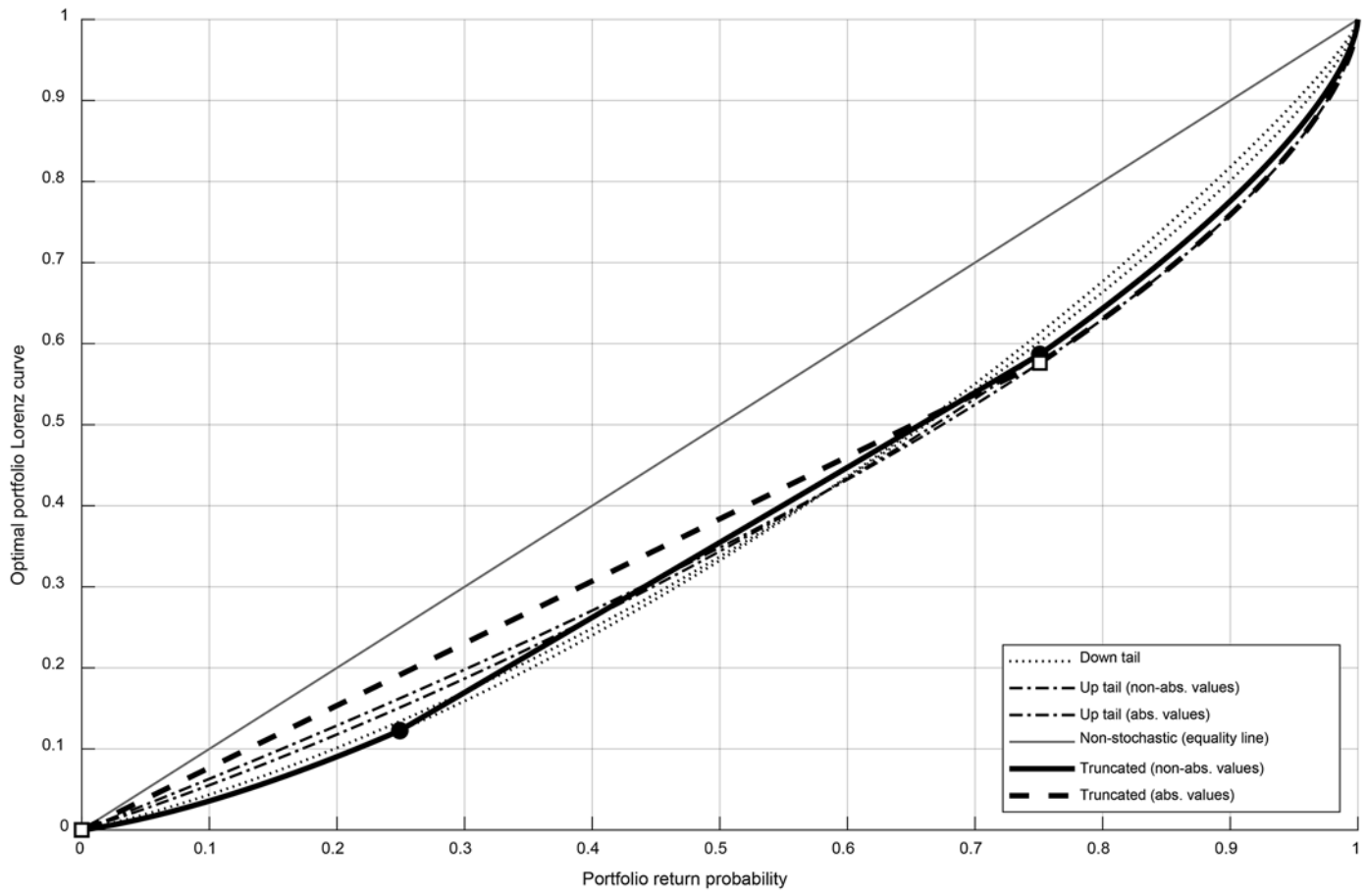

**Fig. 5**. Optimal portfolio's Lorenz curve

The first Lorenz curve is depicted in continuous tick black color. Concretely, in dotted thin black line, a weighted average of the power-law distribution and the Kumaraswamy distribution is shown with some illustrative weights of 0.7 and 0.3 respectively. Up to a sample cumulative return probability of 0.25, the curve characterizes the desired down tail of the optimal targeted Lorenz curve.

---

ature (e.g. *Chen, Tsaur* & *Rhai*, 1982). In an ALC form (Absolute Lorenz curve), i.e., without scaling by the mean, it is considered in (*Shalit* & *Yitzhaki*, 1984). Even in its generalized form, it must be emphasized that the convergence results of *Theorem 1* and *Theorem 3* still hold under suitable restrictions. The latter boil down to imposing proper truncation on the generalized Lorenz curves so that they would be suitable for our iteration procedure. This means that both proofs must remain technically valid, and the modified curves maintain economic sense. Good candidates would be either to take only the increasing section of the generalized Lorenz curve when forming its inverse, or more stringently, its positive increasing part. This allows to take the inverses and also to view the curve as a d.f. In both procedures, after a certain iteration, the Lorenz curves become monotonic and non-negative, and thus classical in shape. Certainly, the economic sense of the limiting distribution, even when reached after an iteration start based on generalized curves, is preserved. Namely, it manifests a special distribution. It is non-negative, as the returns of an optimal portfolio should be, and has some desirable stability properties.

[19] A useful reference for the geometry of Lorenz curves is (*Kämpke* & *Radermacher* (2015), Sections 3.1-3.2, pp. 29 – 37).

Analogously, in dash-dotted thin black line, a weighted average of the power-law distribution and the Kumaraswamy distribution with some illustrative weights of 0.2 and 0.8 respectively is plotted. From sample cumulative return probability of 0.75 to 1, the curve characterizes the desired upper tail of the optimal targeted Lorenz curve. The belly between 0.25 and 0.75 cumulative return probability is a straight line characterizing the optimal targeted Lorenz curve in that segment. The optimal curve is a union of the three segments and is given in continuous tick black color. In terms of the parent distribution, the setting is effectively: (i) a targeted down tail with higher weight on the power-law than the Pareto one; (ii) a targeted upper tail with higher weight on the Pareto law than the power one; (iii) targeted non-stochastic belly in line with SSD[20]. The aim is the Lorenz curve of the optimized portfolio to be as close as possible to the optimal first Lorenz curve from Fig. 5.

The second Lorenz curve is depicted in dashed tick black color. Concretely, it is based on no special down tail and fully on the Kumaraswamy distribution as an upper tail (respective weight of 1). The belly encompasses the down tail and is represented as a straight line from the origin to a point representing a sample cumulative return probability of 0.75 after which the assigned upper tail starts. The auxiliary dotted curves are the same as before.The optimal curve is a union of the two segments and is given in dashed tick black color. In terms of the parent distribution, the setting is effectively: (i) no distinguishable down tail; (ii) a targeted upper tail with higher weight on the Pareto law than the power one; (iii) targeted non-stochastic belly in line with SSD. This curve can be viewed as a Lorenz curve characterizing the absolute returns of the optimal portfolio. It provides an appropriate both theoretically and empirically justified way to avoid the inconvenience to work with generalized Lorenz curves of distributions that can take negative values. The target is again the Lorenz curve of the optimized portfolio to be as close as possible to the optimal Lorenz curve from Fig. 5.

In addition to the above technical motivation driven by the pure Lorenz curve's construct restrictions, some further supporting arguments can be given. It should be noted that often the empirical literature on estimation of asset returns distributions relies exactly on absolute values. This is usually done to get more robust estimates under general L1 norm considerations. In addition, there could also be a special focus on the tail characteristics of the returns. Under such, absolute values work well for fatness classification. This is exactly also the approach of the popular review study of (*Gabaix*, *Gopikrishnan, Plerou* & *Stanley*, 2003) where such are used. Therefore, our measure construct approach is in line and compatible. Furthermore, in the article, Paretian tails are shown to be the norm for various asset returns and that was communicated above to be supporting evidence for the use of our iterated Lorenz curves results from the body text. However, it must be noted that not only does the discovered power law matters but also its precise exponent. Close to cubic ones are detected. The one plus golden section exponent of ours from *Theorem 3* is in line and this holds exactly under an absolute values of returns measurement. All of this adds further structure to the proposed risk measure. Interestingly, various phenomena observed in both social and natural sciences, which share heavy-tailed distributions, appear to have much more in common and thus be tied together. The underlying forces driving them give rise to clear statistical patterns. This enables these phenomena to often be formulated as laws, principles, and rules – such as Zipf's law, Zipf-Mandelbrot's law, and the 80/20 principle, among others. As the foreword of (*Arnold*, 2012) points out, the Pareto-like distributions are a key modeling tool for all of them.

[20] Technically, a straight-line Lorenz curves characterizes flat segments in the parent cumulative distribution function. Relevant examples that could serve as exercises can be found in (*Kämpke* & *Radermacher*, 2015).

To be consistent with the GMD case[21], where essentially the distances are measured by an $l_1$ Kantorovich probability metric (discussed in (*Rachev, Stoyanov* & *Fabozzi*, 2008) in applied context and in (*Rachev, Klebanov, Stoyanov* & *Fabozzi*, 2013) in encyclopedic theoretical one), we will proceed in the same way with the two Lorenz curves from above. This means that up to scaling in Fig. 5, instead of the area between portfolio Lorenz curve and the diagonal in continuous thin black color, we are considering the area between the same portfolio Lorenz curve and any of the prior elaborated optimal ones in tick color. We must emphasize that unlike the GMD, which as discussed at length in (*Furman, Wang* & *Zitikis*, 2017), is a coherent risk measure, we cannot say much a priori about the properties of our new risk measure. It seems to resemble, on one hand, the tail-Gini risk measures (details again in (*Furman, Wang* & *Zitikis*, 2017)), which similarly to our case, are based on Lorenz curve tail integration, but unlike it, the underlying distance is not taken under modulus since the diagonal in Fig. 5 always dominates any curve. On the other hand, there is also resemblance to the well-known in the literature benchmark tracking problems (details in *Rachev, Stoyanov* & *Fabozzi*, 2008) and deviation risk measures (e.g. *Rachev, Stoyanov* & *Fabozzi*, 2008) and (*Rockafellar, Uryasev* & *Zabarankin*, 2002). Last but not least, our measure can also be viewed as a composition one in the spirit of (*Righi*, 2018) and (*Righi*, 2019) since it is a sum of three elements – the belly and the two tail segments. We won't explore comprehensively in this application sup-port appendix the properties of our measure in terms of coherence and fit to the mentioned well-known other risk measures in the literature. For our portfolio optimization purposes, it is sufficient to provide the intuition behind the measure, relate it to the Lorenz curves explored, and present some implementation details. We have done the former two, now we turn attention to the latter.

[21] Other probability metrics are also possible and as we will see are even desirable. Yet, for comparison purposes, we stick to the same metric the GMD in its integral form uses. A good plausible choice is also the Levy metric since as discussed in (*Chobanov* & *Ignatov*, 2002), (*Rachev, Klebanov, Stoyanov* & *Fabozzi*, 2013), and (*Rachev, Stoyanov* & *Fabozzi*, 2008), it metrizes the weak convergence, but the problem with it is its complex technical treatment. The weak convergence can be useful for more flexible approximations.

### 3. Technical set-up

We continue with the concrete mathematical construct of our new risk measure. We will name it a GS measure ($GS_1$ – when the first optimal Lorenz curve is used, and $GS_2$ when the second one). As noted, it relies on the areas from Fig. 5 previously discussed which can easily be computed. We would like to impose equal scaling between the GS measure and the GMD one. After this introduction, we move to the necessary notation. Let $N$ be the number of assets in a portfolio $a$ under focus and denote by $r = [r_1, ..., r_N]$ the vector of their returns over a fixed horizon. The portfolio's return is $r_\alpha = w_1r_1 + w_2r_2 + ... + w_nr_n$, where $w = [w_1, w_2, ..., w_n]$ is the vector of weights on the individual assets. The GMD measure has the form:

$$GMD(r_\alpha) = E\left|r_\alpha^{(1)} - r_\alpha^{(2)}\right|, \tag{77}$$

where $r_\alpha^{(1)}$ and $r_\alpha^{(2)}$ are two i.i.d. sample returns from $r_\alpha$. It is well known (e.g. *Arnold*, 2015; *Furman, Wang* & *Zitikis*, 2017; *Mansini, Ogryczak* & *Speranza*, 2007; and *Shalit* & *Yitzhaki*, 1984) that if by $G(r_\alpha)$ we denote the Gini's index of $r_\alpha$, then the following holds:

$$G(r_\alpha) = \frac{GMD(r_\alpha)}{2E(r_\alpha)} = \frac{E\left|r_\alpha^{(1)} - r_\alpha^{(2)}\right|}{2\mu_{r_\alpha}} = 2\int_0^1 \left[x - L_{r_\alpha}(x)\right]dx, \tag{78}$$

where $\mu_{r_\alpha} = E(r_\alpha)$. Now, if we denote by $S$ the area between the curve in (continuous/dashed) tick black color and the diagonal in continuous thin black color from Fig. 5, we will get:

$$\begin{aligned} S &= \int_0^{\beta_{down}} \left[x - L_{\text{down}}(x)\right]dx + \int_{\beta_{\text{down}}}^{\beta_{up}} \left[x - L_{\text{center}}(x)\right]dx + \int_{\beta_{up}}^{1} \left[x - L_{up}(x)\right]dx \\ &= \int_0^1 \left[x - L_{GS}^{\beta,\gamma}(x)\right] \end{aligned} \tag{79}$$

where $\beta_{down}$ and $\beta_{up}$ are the cut-off probabilities defining the tails of the optimal portfolio's Lorenz curve ($\beta = [\beta_{down}, \beta_{up}]$).

If we set by $\gamma_{dTailPa}$ and $\gamma_{dTailP}$ the down tail weights of the Pareto distribution and the power-law distribution induced Lorenz curves respectively, with $\gamma_{uTailPa}$ and $\gamma_{uTailP}$ being the counterparts of the latter two weights for the upper tail ($\gamma = [\gamma_{dTailPa}, \gamma_{dTailP}, \gamma_{uTailPa}, \gamma_{uTailP}]$), we have additionally:

$$L_{\text{down}}(x) = \gamma_{d\text{TailPa}}\left(1-(1-x)^{\frac{\sqrt{5}-1}{2}}\right) + \gamma_{d\text{TailP}}\, x^{\frac{\sqrt{5}+1}{2}} \tag{80}$$

$$L_{up}(x) = \gamma_{uTailPa}\left(1-(1-x)^{\frac{\sqrt{5}-1}{2}}\right) + \gamma_{uTailP}\, x^{\frac{\sqrt{5}+1}{2}} \tag{81}$$

$$L_{\text{center}}(x) = \frac{L_{up}\left(\beta_{up}\right) - L_{\text{down}}\left(\beta_{\text{down}}\right)}{\beta_{up} - \beta_{\text{down}}}\left(x - \beta_{\text{up}}\right) + L_{up}\left(\beta_{up}\right) \tag{82}$$

$$\begin{aligned} L_{GS}^{\beta,\gamma}(x) = {} & L_{\text{down}}(x) 1_{\{x\in[0,\, \beta_{\text{down}})\}} + L_{\text{center}}(x) 1_{\{x\in[\beta_{\text{down}},\, \beta_{up})\}} + \\ & + L_{up}(x) 1_{\{x\in[\beta_{up},\, 1]\}} \end{aligned} \tag{83}$$

Our $GS_1$ risk measure is defined by:

$$GS_1\left(r_\alpha\right) = \frac{2\mu_{r_\alpha}}{\int_0^1 L_{GS_1}^{\beta_1,\gamma_1}(x)dx} \int_0^1 \left| L_{r_\alpha}(x) - L_{GS_1}^{\beta_1,\gamma_1}(x) \right| dx, \tag{84}$$

where we put the superscripts of β and γ to differentiate the measure despite the no restrictions imposed on the coefficients and additionally we clearly have $L_{GS_1}^{\beta_1,\gamma_1}(x) = L_{GS}^{\beta,\gamma}(x)$.

For the $GS_2$ case, imposing the restrictions $\beta_{down} = 0$, $\gamma_{uTailP} = 1$, and $\gamma_{uTailP} = 0$ in (80), (81), (82), and (83), we get:

$$L_{GS_2}^{\beta_2,\gamma_2}(x) = L_{\text{center}}(x) 1_{\{x\in[0,\beta_{2,up})\}} + L_{up}(x) 1_{\{x\in[\beta_{2,up},1]\}} \tag{85}$$

and thus also:

$$GS_2\left(r_\alpha\right) = \frac{2\mu_{|r_\alpha|}}{\int_0^1 L_{GS_2}^{\beta_2,\gamma_2}(x)dx} \int_0^1 \left| L_{|r_\alpha|}(x) - L_{GS_2}^{\beta_2,\gamma_2}(x) \right| dx, \tag{86}$$

where we made also explicit the use of absolute value of the portfolio returns. In both cases, inside the integral, we have the targeted area from Fig. 5. The front coefficient is a scaler guaranteeing consistency with the GMD.

It should be emphasized that for both risk measures, we again have representations similar to those in (75) and (76). Specifically, formulas (84) and (85) allow both measures to be interpreted as weightings on risk scenarios, though the exact weights have more complex shapes.

Although Gini's index was presented in its classical form in (78), which helped define our risk measures $GS_1$ and $GS_2$, it is important to note that, as discussed in (*Aaberge*, 2000) and (*Arnold*, 2015), significant generalizations of the index exist. A popular such involves incorporating risk aversion which allows for greater flexibility in calibrating investor risk preferences. As shown in the next section, this is crucial for producing realistic empirical results. Risk aversion can be introduced by appropriately weighting the risk scenarios in addition to the one coming from the intrinsic differences among the measures considered. Such parameterized scaling can be left to the modeler's discretion, aside from the selection of the risk measure. A possible approach was introduced in (*Donaldson* & *Weymark*, 1980), leading to a variety of applications as summarized in (*Aaberge*, 2000) and (*Arnold*, 2015). For example, the technique was applied in portfolio theory in (*Shalit* & *Yitzhaki*, 1989) and (*Yitzhaki*, 1983), building on the foundational work in (*Shalit* & *Yitzhaki* (1984), where the GMD methodology in finance was first introduced.

Following the exposition of (*Aaberge*, 2000), extended GMD and Gini's indices can be composed in the form:

$$G(r_\alpha, v) = \frac{GMD(r_\alpha, v)}{2E(r_\alpha)} = v(v-1)\int_0^1 (1-x)^{v-2}\left[x - L_{r_\alpha}(x)\right]dx, \tag{87}$$

where the weight $v(v-1)(1-x)^{v-2}$ is associated with each section of the area between $x$ and $L_{r_\alpha}(x)$. As discussed in (*Aaberge*, 2000) and (*Arnold*, 2015), they can be viewed as linear combinations of moments generated by the Lorenz curve when interpreted as a distribution function[22]. This allows the extended in-dex $G(r_\alpha, v)$ to be considered a good instrument to capture risk aversion. When $v = 2$, we obtain the standard Gini index. Additionally, for $v = 1$, we have the case of a risk-neutral investor. As $v \to +\infty$, the extreme case of a max-min investor emerges. Finally, for $0 < v < 1$, we encounter the opposite extreme – a risk-loving investor.

[22] More precisely, $G(r_\alpha, \upsilon)$ is a linear combination of the terms $D_v(r_\alpha) = \frac{v+1}{v}\int_0^1 x^v dL_{r_\alpha}(x) - \frac{1}{v}$, where $\int_0^1 x^v dL_{r_\alpha}(x)$, known in the literature as LC moments (see *Aaberge*, 2000), are appropriately scaled to achieve uniform range. The exact form of the representation is: $G(r_\alpha, v+1) = 1 + (v+1)\sum_{i=1}^{v}(-1)^i \binom{v}{i}\frac{i}{i+1}(1 - D_i(r_\alpha))$. The family $\{D_v(r_\alpha)$: $v = 1,2,...\}$ can itself be viewed as a set of risk measures. Since the curve $L_{r_\alpha}(x)$, in its d.f. form, is defined on a bounded interval, it is uniquely determined by its moments. Thus, the measures $D_v(r_\alpha)$ fully determine $L_{r_\alpha}(x)$. In its Lorenz curve form, the latter, along with the mean $E(r_\alpha)$, completely determines the distribution of $r_\alpha$.

The above exposition allows us to formulate the extended $GS_1(r_\alpha, v)$ and $GS_2(r_\alpha, v)$ measures accordingly:

$$GS_1\left(r_\alpha, v\right) = \frac{\mu_{r_\alpha} v(v-1)}{\int_0^1 L_{GS_1}^{\beta_1,\gamma_1}(x)dx} \int_0^1 (1-x)^{v-2} \left| L_{r_\alpha}(x) - L_{GS_1}^{\beta_1,\gamma_1}(x) \right| dx \tag{88}$$

$$GS_2\left(r_\alpha, v\right) = \frac{\mu_{|r_\alpha|} v(v-1)}{\int_0^1 L_{GS_2}^{\beta_2,\gamma_2}(x)dx} \int_0^1 (1-x)^{v-2} \left| L_{|r_\alpha|}(x) - L_{GS_2}^{\beta_2,\gamma_2}(x) \right| dx \tag{89}$$

They are the ones we will mainly employ in our consequent empirical analysis.

### 4. Model implementation

We focus on solving classical portfolio optimization problems as in (*Rachev, Stoyanov & Fabozzi*, 2008). Namely, for a given expected return, minimization of a risk measure is implemented. This allows to get the optimal risk-return profiles, the efficient frontiers, and the optimal weights. We use scenarios based on either directly historical returns as in (*Cesarone*, 2020), or first estimation of suitable copulas and then generation of samples as in (MathWorks Quant Team, 2024). Additionally, we implement optimizations considering several popular risk measures: variance, conditional value at risk (CVAR), mean absolute deviation (MAD), Gini's mean difference (GMD), and finally our golden section ones ($GS_1$) and ($GS_2$). This is done for comparison purposes.

The books of (*Cesarone*, 2020) and (*Mansini, Ogryczak, Speranza* & EURO, 2015) are good references on numerical portfolio optimization. The latter has special emphasis on formal reformulation of various risk minimization problems into the language of linear programming (LP). This is very important for getting more runtime efficient solutions based on variations of the simplex method (primal or dual). The former gives a more modern perspective of the same together with implementations in Matlab. As already indicated, with the practical example of the appendix, we aim mainly illustration and not delving into specialized portfolio theory agenda. So, we won't focus on the most efficient numerical procedures which can be suitable for such. Our online appendix and the companion code provide programmatic realizations that are both highly operative and parsi-monious and can easily be followed[23]. There are enough to back our analysis and conclusions.

Let's turn attention now to the exact mathematical form of our optimization problem. For the $GS_1$ measure, it is:

$$\mathrm{Min}_{\mathrm{W}} \frac{1}{T} \frac{2\mu_{r_\alpha} v(v-1)}{\int_0^1 L_{GS_1}^{\beta_1,\gamma_1}(x)dx} \sum_{i=1}^{T} \left(1-\frac{i}{T}\right)^{v-2} \left| \frac{\sum_{j=1}^{i} r_\alpha^{(j:T)}}{\sum_{j=1}^{T} r_\alpha^{(j:T)}} - L_{GS_1}^{\beta_1,\gamma_1}\left(\frac{i}{T}\right) \right| \tag{90}$$

$$\text{s.t. } \mu_{r_\alpha} = \sum_{i=1}^{n} w_i \mu_i \tag{91}$$

$$\sum_{i=1}^{n} w_i = 1, \tag{92}$$

where for fixed weights $w$, $r_\alpha^{(j:T)}$ are ordered sample returns from $T$ scenarios. We can assume that each of the latter has a $\frac{1}{T}$ probability of occurrence. Effectively, we discretize the two Lorenz curves participating in the measure and use that the optimal one has a predetermined mathematical form. Analogously, for the $GS_2$ measure, the problem is:

$$\mathrm{Min}_{\mathrm{W}} \frac{1}{T} \frac{2\mu_{|r_\alpha|} v(v-1)}{\int_0^1 L_{GS_2}^{\beta_2,\gamma_2}(x)dx} \sum_{i=1}^{T} \left(1-\frac{i}{T}\right)^{v-2} \left| \frac{\sum_{j=1}^{i} |r_\alpha|^{(j:T)}}{\sum_{j=1}^{T} |r_\alpha|^{(j:T)}} - L_{GS_2}^{\beta_2,\gamma_2}\left(\frac{i}{T}\right) \right| \tag{93}$$

$$\text{s.t. } \mu_{r_\alpha} = \sum_{i=1}^{n} w_i \mu_i \tag{94}$$

$$\sum_{i=1}^{n} w_i = 1, \tag{95}$$

with the difference that now $|r_\alpha|^{(j:T)}$ are ordered sample returns in absolute value from $T$ scenarios.

The optimization in the above two problems is non-linear. It is such due to not only all the absolute values that participate in the objective function but also to the ordered statistics sums presence since ordering is a non-linear operation. Following similar logic to the one from (*Cesarone*, 2020) and (*Mansini, Ogry-*

[23] We base our empirical results on them. The code has also other realizations for runtime comparisons as supplementary material. They might be continuously improved following a more specialized portfolio theory research path. For the main realizations, we use the Matlab's build-in non-linear solvers and the free CVX package.

*czak*, *Speranza* & EURO, 2015) we can linearize the above problems[24]. This would be done at the expense of notation parsimony and a big number of helper variables, with the latter certainly being a drawback to the run-time, but still keeping it to acceptable limits. As already noted, we will prefer for this research more straightforward solutions. Namely, since programmatically the above problems can easily be coded, we will directly achieve the posed parsimony goal. That would be at the expense of the non-linear optimization routines we must resort to and the higher expected runtime than the one the LP might offer. We experimented with different solvers, and the one chosen – despite generally not being recommended for non-smooth functions – provides decent results, as we will see in the next section.

Furthermore, we might modify the above basic formulations of the optimization portfolio problem by adding additional constraints. Namely, we always require positive values for $\mu_{r_a}$ since our employed Lorenz curve definition relies on that. Additionally, we might also pose to work with positive weights to better explore the winner-takes-all problem.

Finally, we need to discuss the sequence of the optimizations done. At the first step, for all the problems considered, we run an optimization without posing a prespecified expected return as an equality constraint. This allows to find from the optimal weights both the expected return (weighted average of the expected return of the individual assets) and the risk measure (value of the objective function). The latter two determine the starting point of the efficient frontier. At the second step, we determine the expected returns vector under interest. We take as its minimum the just found expected return from the first step, and as its maximum, the expected return of the asset with the highest such. Then we discretize the interval between the minimum and the maximum with the number of points we would like to have at the efficient frontier. At the third step, we run a sequence of constraint optimizations by posing expected returns equal to each entry of the determined at the previous step vector of expected returns as an equality constraint.

## 5. Empirical results

Our data is based on the constituents of two popular international stock indices. We consider the Dow Jones Industrial Average (DJI) and the EURO STOXX 50. The former is composed of 30 stocks of large US companies which are also blue-chips and are considered representative for the US economy. The latter is composed of 50 stocks of some of the largest companies of 11 eurozone countries. Again, these are blue-chips and can be viewed as representative for the euro area. For each of the two indices, we select a subsample from their respective constituents. Then we try to build USD and an EUR denominated optimal portfolios respectively. Our scenario generation procedure would require the use of as long as possible historical data series. Only this would allow a precise estimation of the joint empirical densities needed. Certainly, we cannot consider all the constituents of each index for our portfolios. First, this is unnecessary for our analysis to draw valid conclusions and would only lead to an inefficient use of computational resources. Second, more prosaically, the data itself does not allow comprehensive inclusion. The indices update their entries, and this does not make possible to have simultaneously long time series for each stock. Data cleaning is needed and the companies with insufficient history have to be excluded.

[24] As well as the ones elaborated to be considered: CVAR, MAD, and GMD. For mean-variance, the norm is the use of quadratic programming although non-standardized LP optimizers also exist.

Concretely, we use data from LSEG's Refinitiv Workspace (formerly Eikon)[25]. We implement the data cleaning in several steps. First, we may observe from the provided data files that as expected, each index is subject to an active process of companies joining and leaving. We can get all the unique stocks that have ever been in the index and extract data for them. This is our step 0. Then in step 1, we drop all the companies that have not been long enough in the index. For DJI, the data provider distributes quality data since Jan 1, 1994. We take series from that date till June 30, 2024. We may pose to work with stocks that cover the period almost exclusively, say, more than 95% of time. We approximate that span to 7670 trading days. After imposing this constraint, the 56 starting stock entries reduce to 40. Finally, in step 2, we check whether all the stocks selected cover the same trading days with also having liquid quotes for each time point. If there is no quote for specific date for some stock, we remove that date from the whole sample. This produces a final sample of 40 stocks over 7650 trading days.

For EURO STOXX 50, we have data from the data vendor since August 12, 1998, and again take series till June 30, 2024. We do that for all the unique stocks that have ever been in the index. We follow a similar type of data cleaning procedure as for DJI. We drop all the stocks that have less presence than 95% of the time. We approximate that to 6435 days. Since the index experienced very active updates, from the 124 unique entries, we get only 42 that fulfil the requirement. Additionally, we pose the condition to also have coverage of the most recent years because they are formative for our scenarios. Finally, we again drop the trading days that are not common for all the selected stocks. Our final sample is of 42 stocks over 6098 trading days.

The risk minimization procedures are implemented in Matlab. As elaborated before, we focus on variance, CVAR, MAD, GMD, $GS_1$ and $GS_2$ risk measures.

[25] Stocks data is easily obtainable from a variety of providers without charge.

The Financial Toolbox has ready object-oriented implementations of the first three based on classical papers (see the Mathworks documentation[26]). Possibly following the lack of consensus in the literature for the best optimization method for GMD, Matlab does not provide implementation[27] of that risk measure. We devise such based on our $GS_1$ methodology[28] from *Section 3* just by taking $L_{GS}^{\beta,\gamma}(x) = x$ in (87) and also noting that $\int_0^1 L_{GS_1}^{\beta_1,\gamma_1}(x)dx = 0.5$. Both for the $GS_1$ and $GS_2$ risk measures, we do non-linear optimization. Similarly to (*Rachev, Stoyanov & Fabozzi*, 2008, Section 9.5, p. 303 and Section 9.7.5, p. 314), at the first step, we pose equal weights as initial guess. For all risk measures, we consider 10 portfolios for the efficient frontier as in (MathWorks Quant Team, 2024) and (*Rachev, Stoyanov & Fabozzi*, 2008, Section 9.5, p. 301) which are enough to get appropriate visualizations and draw relevant conclusions. This leads to 10 consecutive optimizations. For each step, we will use as initial weights guess the obtained solution from the previous one as in (*Rachev, Stoyanov & Fabozzi*, 2008).

Each measure, as previously indicated, captures different aspects of investor preferences uniquely. This is reflected in distinct vectors of weights among the optimal portfolios. They would be convenient to see not only technically how the optimization worked and whether we got reasonable values but also a set of other characteristics.

The first feature is to what extent the optimal portfolios are close to each other both from statistical and financial points of view. It may happen that despite being different, the chosen risk measures perform in a similar fashion because the exact risks they are sensitive to have the same dominant presence in the data. Thus, even being proxies for specific type of investor's preferences, since the quantification of the risks may not be precise and comprehensive by these risk measures, this may trigger similar values of them when in general should not. Vice versa, the measures can be oversensitive to small disparities in the inherent risks of the data. Both situations do not mean necessarily that the chosen risk measures are bad per se. We might just have a wrong selection of them for the environment under scope and an alternative one to have been better. For example, data primarily driven by variability risks will produce similar weights for portfolios optimized under MAD and GMD. Upon that, if we additionally have occurrence of multivariate normality, or even be close to it, both the CVAR and the variance optimal portfolios will be similar to each other as well as to the entities from the previous two cases[29]. In situation of pure tail risks, the variance optimization is completely inadequate and will show results in dissonance with the other discussed risk measures. Among the latter, MAD and GMD will certainly produce distinct output but interestingly not extremely so. The reason is that the two measures succeed to capture both variability and tail risks, yet, in different ways. The CVAR optimal weights may or not be much different to them. Both the exact composition of the two types of risk matters and the very specifics of each of it. Factors such as the tails variability, tails thickness, distribution skewness, and the disparity from normality kick in among others.

[26] More precisely, the bibliography section at: https://www.mathworks.com/help/finance/bibliography.html#bsw8cho-1.

[27] As of version R2024a – the one used in the paper.

[28] As elaborated before, for comparison reasons, we may also consider other methods, as well as computer languages, for all the risk measures in the extended online appendix coverage. They are not affecting the results of the paper which are based solely on the implementation described in it.

To see the differences between the risk measures according to the logic discussed above, it would be instructive to compare the efficient frontiers they produce. They will give the right visual perspective for closeness and disparity. Additionally, a uniform unit of measurement for risk is essential to enable accurate comparisons. A way-out is to choose a specific risk measure and then within the risk-return coordinates its efficient frontier lies, to represent the efficient frontiers generated by the other measures. This is straightforward to be done. We can just take the 10 vectors of optimal weights forming the frontier of each measure and see what portfolios they generate in the fixed coordinate system of the set benchmark one. Certainly, within the latter, these portfolios will be suboptimal since the optimization criterion was not the one applied to them. Yet, risk-return profiles will be generated represented by dots that can be connected. These new curves will lie below the efficient frontier of the benchmark measure. If there are exceptions, such would be due to numerical issues. The resulting picture gives very good perspective to visually inspect the hinted proximity of the efficient frontiers between the different measures. It should be noted that not only may the usual convexity of the efficient frontier not hold for the mentioned risk curves[30] but their shape may also be irregular in some situations. This is not a problem because the imported objects are just non-optimal portfolios.

The second feature is the exact weights' vector composition of the optimal portfolios and their evolution along the efficient frontier. Some stocks may find marginal to no participation or excessive such. We can plot the weights and see the patterns that emerge. The one under inspection would be the weights' outliers, variability, and stability as well as portfolio concentrations and diversification. Certainly, stability and diversification are desirable properties to be sought. The latter has special sounding if we impose restrictions on the weights. Doing this, effectively reduces the scope for trading which hinders diversification. For example, to better investigate the winner-takes-all problem, we may limit the weights to be positive (i.e., no short selling). Under this setting, for the last points of the efficient frontier, we may get very ostensible concentrations on one highly performing stock due to the less diversification opportunities. So, weight restrictions are not necessarily a bad thing for analysis.

[29] As the example in (MathWorks Quant Team, 2024) essentially shows. Yet, even for that case, very small differences are noticeable.

[30] Yet, we will still call them efficient frontiers hereinafter for convenience.

After the methodological clarifications from the above paragraphs, we move to the presentation of our results. They can be summarized in two figures for each index for each simulation: the composition of the efficient frontier returns across the measures and the efficient frontiers of the measures within the benchmark one.

To produce them, we have to set some further parameters. First, we choose MAD as our benchmark measure. This is due to the fact that it is not Lorenz curve based which prompts a sort of neutrality compared to the alternatives in the list. Additionally, we would expect MAD to perform better than the other measure conceptually close to it in the face of the variance. So, the choice is reasonable. Second, we choose a random selection of stocks from the available set of 40 – we take every third entry. This produces a total number of 14 stocks, which is a consistent span and enough for pertinent analysis. Third, we work with scenarios based on both simulations and historical data. In the former case, for estimation of the empirical copula, we take our whole data sample. Fourth, we consider daily data. Fifth, we stick to 500 scenarios when they are based on historical data. In this case the last observations are taken. We can deal with more scenarios under this setting, but we give a reasonable preference to the latest history. We stick to 1000 scenarios when they are based on simulation. Both numbers are reasonable, and no numerical issues are visible. They allow drawing consistent conclusions under reasonable laptop run-time. Sixth, we use the same weights and cut-off values for the tails as the ones employed in Fig. 5. Seventh, as hinted before, we will work under the restriction of positive weights. This allows a better analysis on the diversification effects. The latter are just more difficult to be traced when there are reinforced by no restraint on short positions. Eighth, for the CVAR, we pose 0.95 which is a reasonable value also used in (MathWorks Quant Team, 2024). Nineth, we use different values for $\upsilon$ to see how the risk aversion affects our analysis.

Our results produce interesting patterns which allow to get some delicate intricacies of the measures under scope. In this line of thought, we can begin with a surprising result. We always find the simple (non-extended) $GS_1(r_\alpha)$ measure to produce very close outcomes to the ones of the simple $GMD(r_\alpha)$ measure, not only along the efficient frontier plots but also in direct comparisons of the numerical values of the weights. The trivial reason is that in all our cases, the Lorenz curves of the span of feasible portfolios always seem to lie below the $L_{GS_1}^{\beta_1,\gamma_1}(x)$ curve and therefore also below the diagonal of the square (equality line) from Fig. 5. Thus, minimizing the simple $GMD(r_\alpha)$, automatically minimizes $GS_1(r_\alpha)$ as well. This is further visible if we look at the Lorenz curves of the optimal portfolios generated by these two measures. With minor differences, the same is valid for the $GMD(r_\alpha, v)$ and $GS_1(r_\alpha, v)$ measures while we maintain an equal value of $\upsilon$ across them for ceteris paribus reasons. Obviously, to be able to better distinguish between the $GMD$ and the $GS_1$ measures, we need to a priori achieve stronger diversification. This can be done only with much larger samples than the 14 stocks considered. Unfortunately, under the currently set methodology of non-linear optimization, this is not feasible with reasonable run-time. Therefore, we can focus on the $GMD(r_\alpha, v)$ and $GS_2(r_\alpha,v)$ measures and see what situation we get there.

Fig. 6 – 7 below display both the compositions of the optimal portfolios and the efficient frontiers for DJI under 500 historical scenarios. Fig. 8 – 9 do the same for that index under 1000 simulation scenarios. Then we proceed to the other index. Fig. 10 – 11 display both the compositions of the optimal portfolios and the efficient frontiers for EURO STOXX 50 under 500 historical scenarios. Fig. 12 – 13 do the same for that index under 1000 simulation scenarios. Everywhere we take $v = 2.5$, which is a plausible value for risk aversion also used in (*Shalit* & *Yitzhaki*, 1989). A sequence of implications follows.

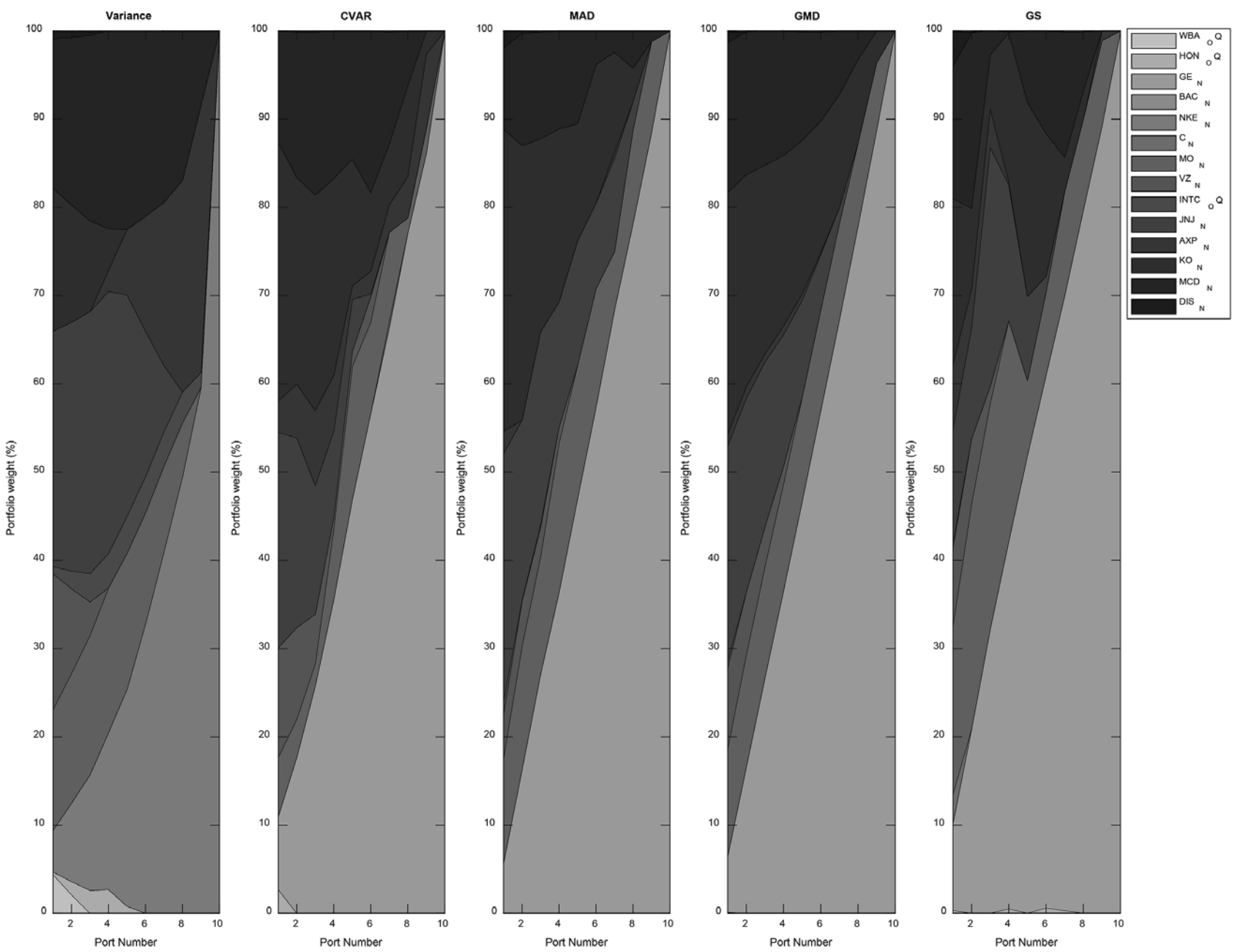

**Fig. 6.** Optimal portfolio compositions DJI, historical scenarios (VAR, CVAR, MAD, GMD, and $GS_2$ )

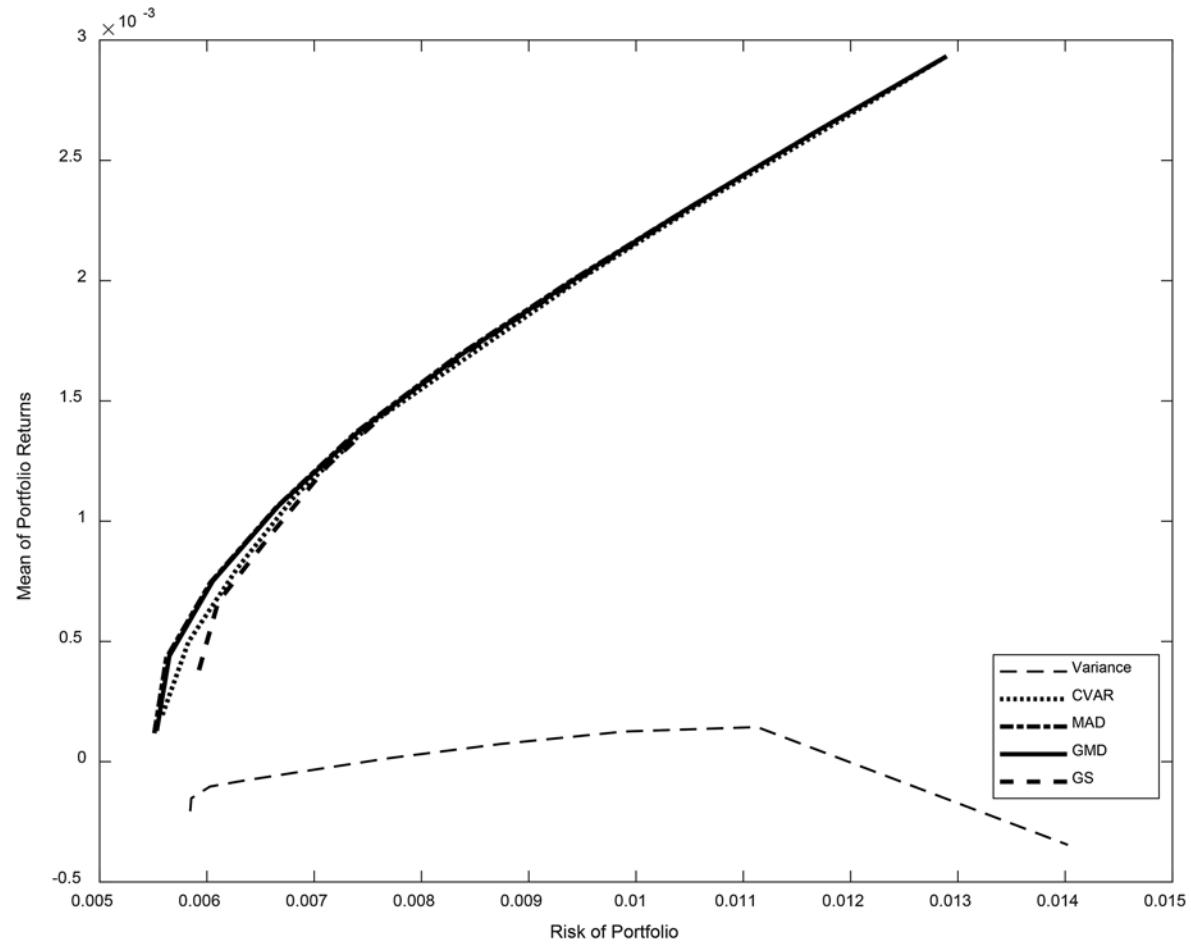

**Fig. 7**. Efficient frontiers DJI, historical scenarios (VAR, CVAR, MAD, GMD, and $GS_2$)

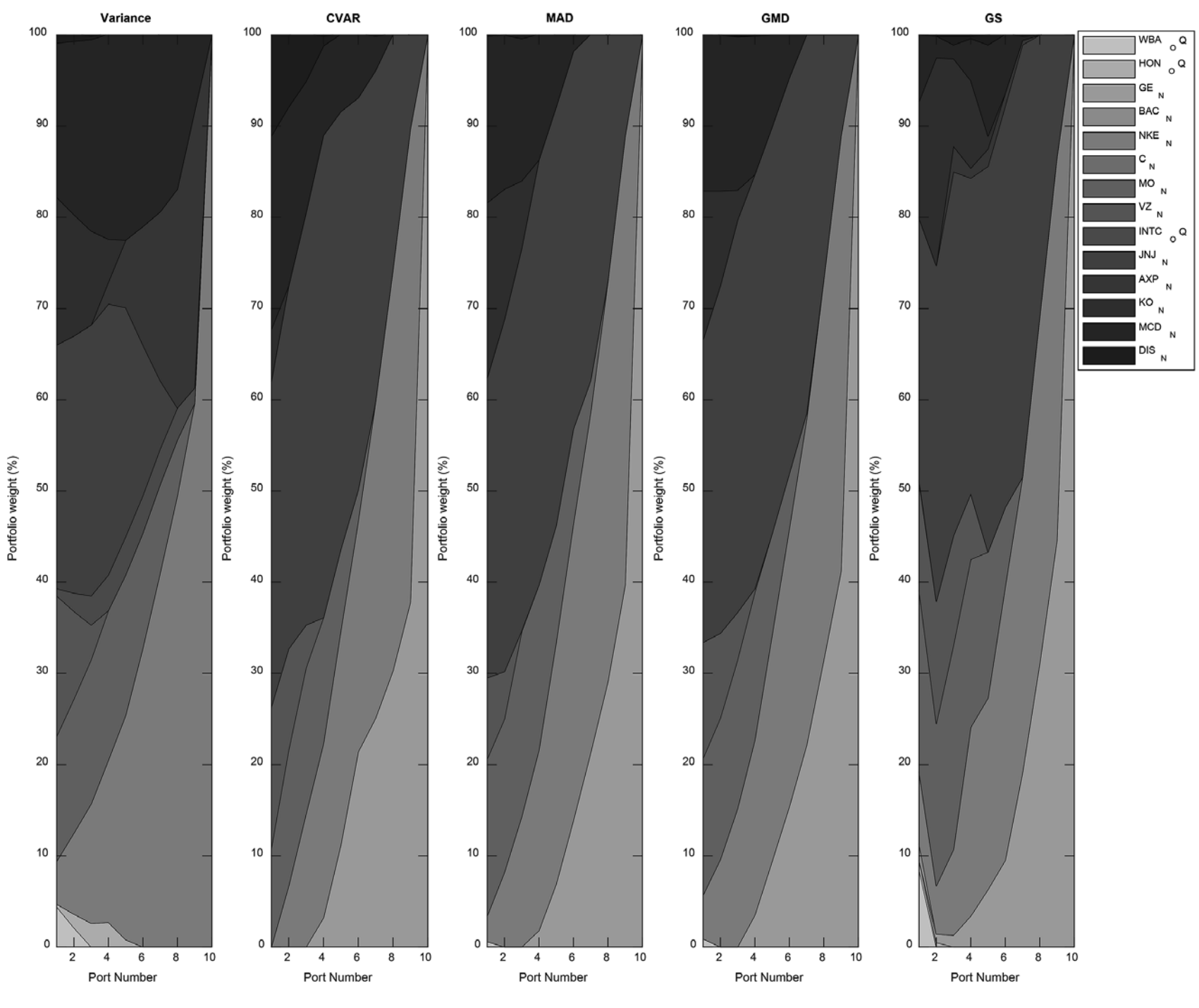

**Fig. 8**. Optimal portfolio compositions DJI, simulation scenarios (VAR, CVAR, MAD, GMD, and $GS_2$)

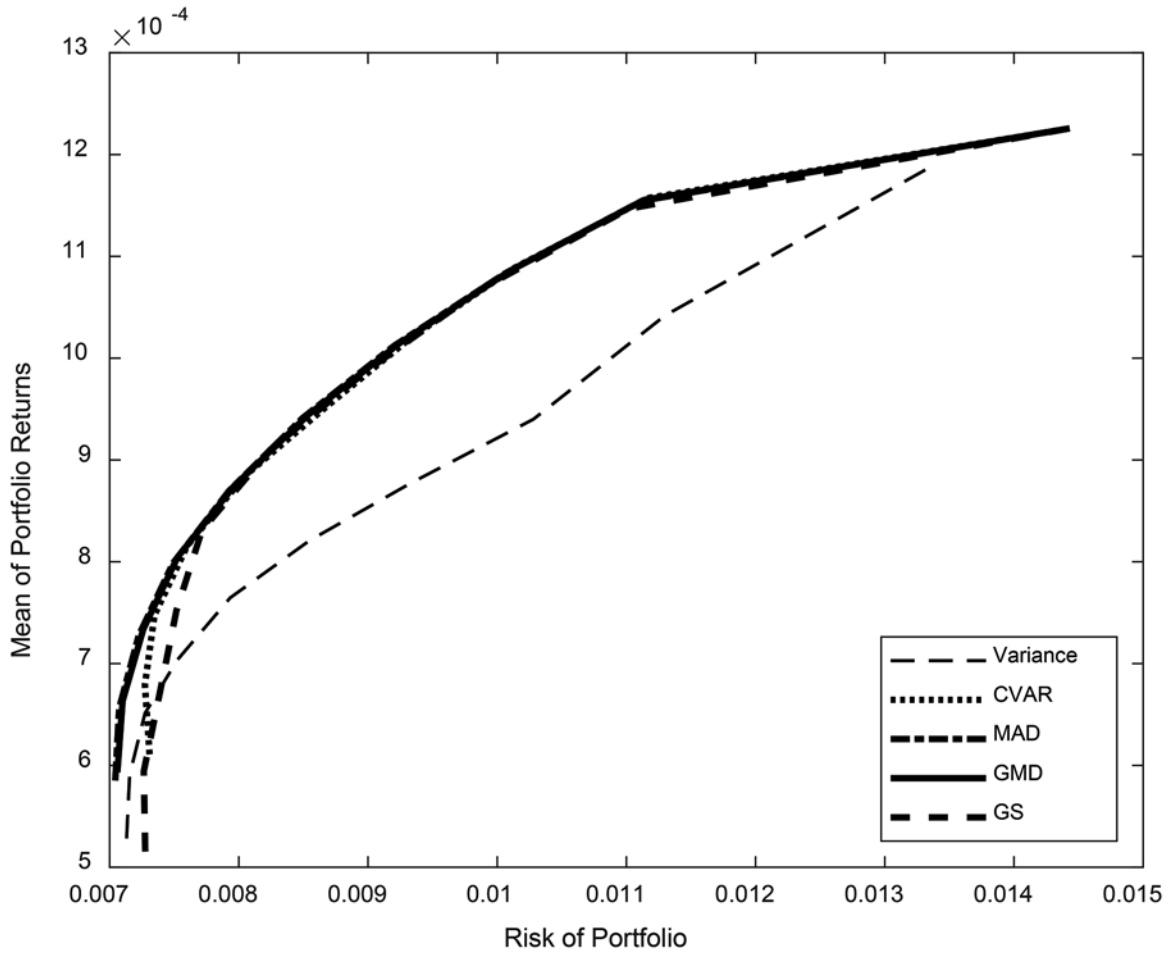

**Fig. 9**. Efficient frontiers DJI, simulation scenarios
(VAR, CVAR, MAD, GMD, and $GS_2$)

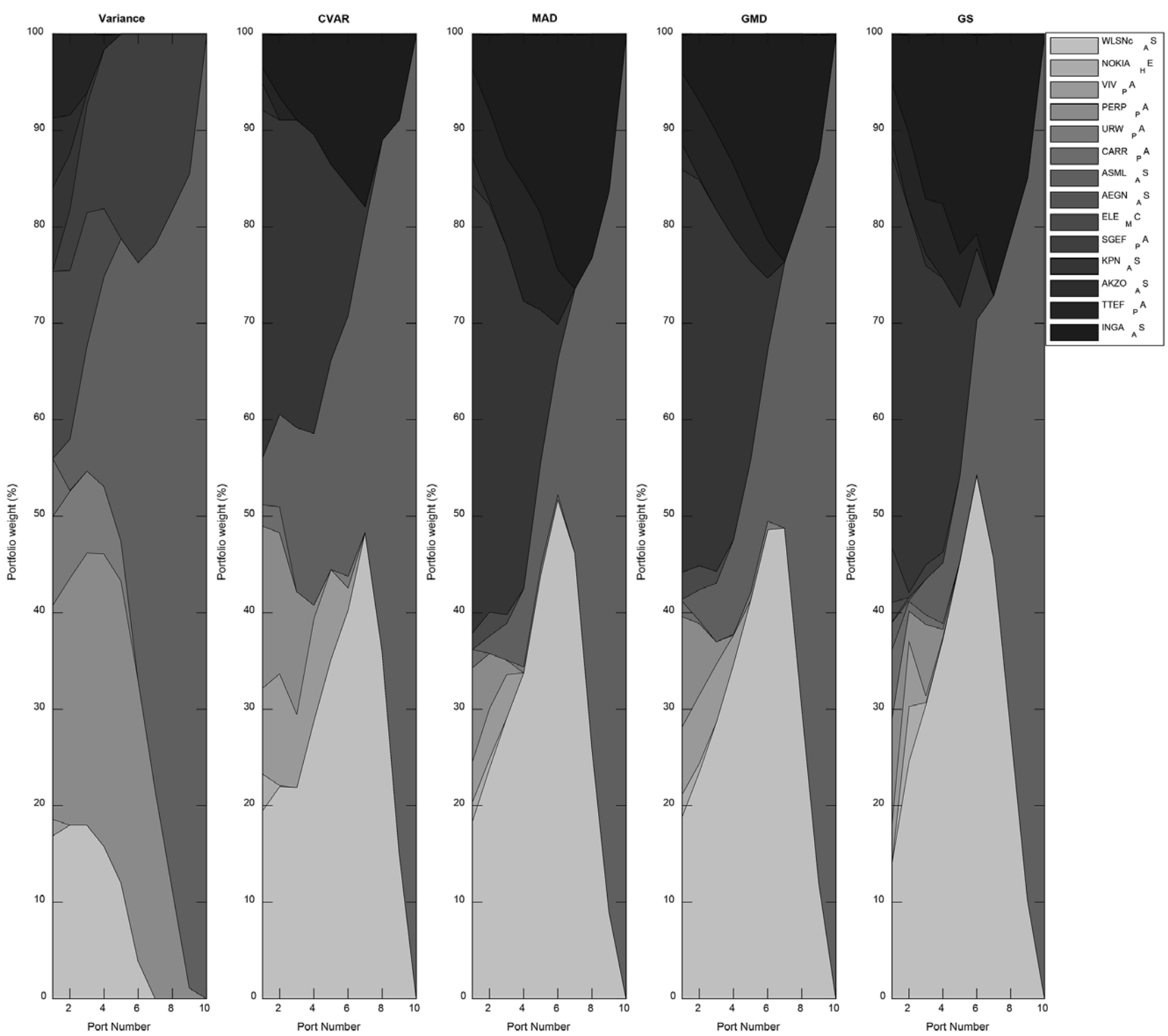

**Fig. 10**. Optimal portfolio compositions EURO STOXX 50, historical scenarios
(VAR, CVAR, MAD, GMD, and $GS_2$)

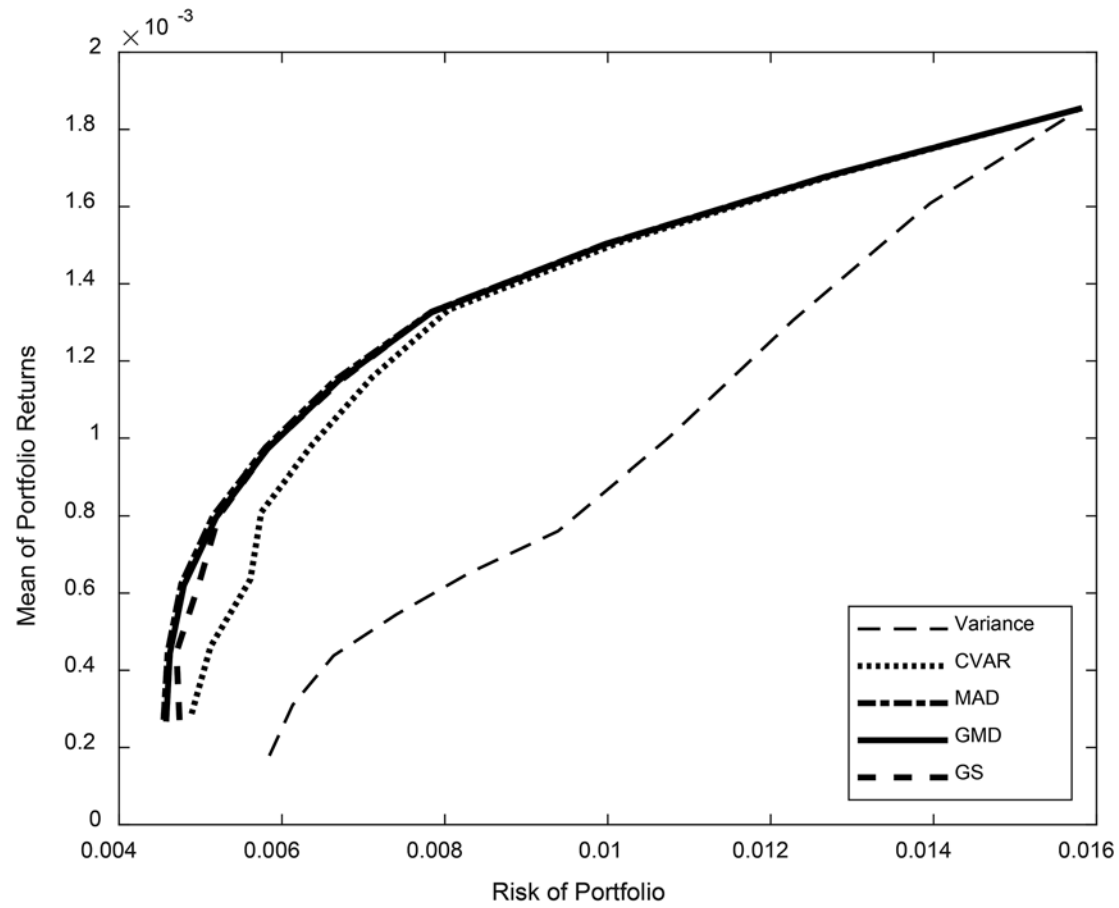

**Fig. 11.** Efficient frontiers EURO STOXX 50, historical scenarios (VAR, CVAR, MAD, GMD, and $GS_2$)

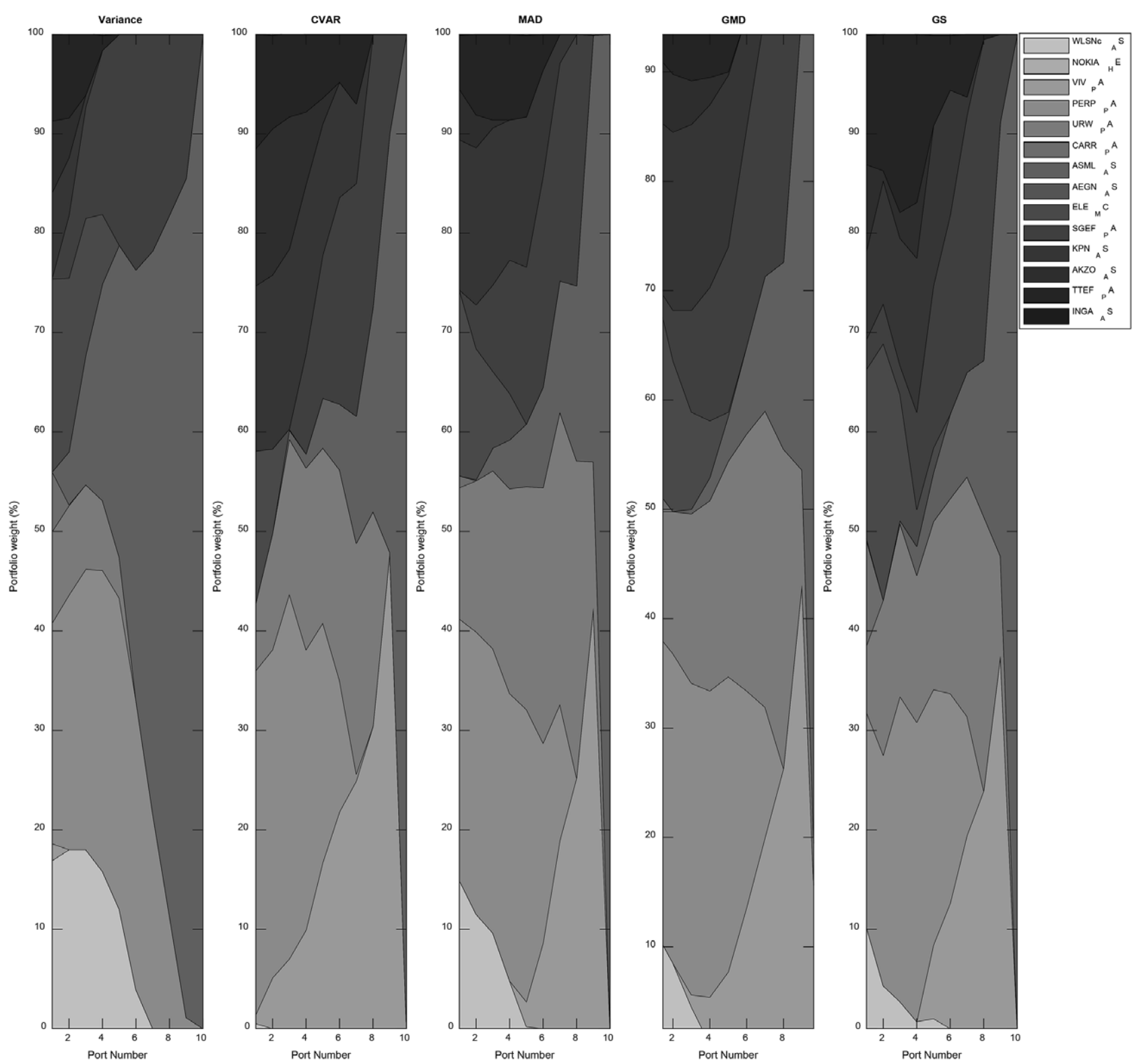

**Fig. 12**. Optimal portfolio compositions EURO STOXX 50, simulation scenarios (VAR, CVAR, MAD, GMD, and $GS_2$)

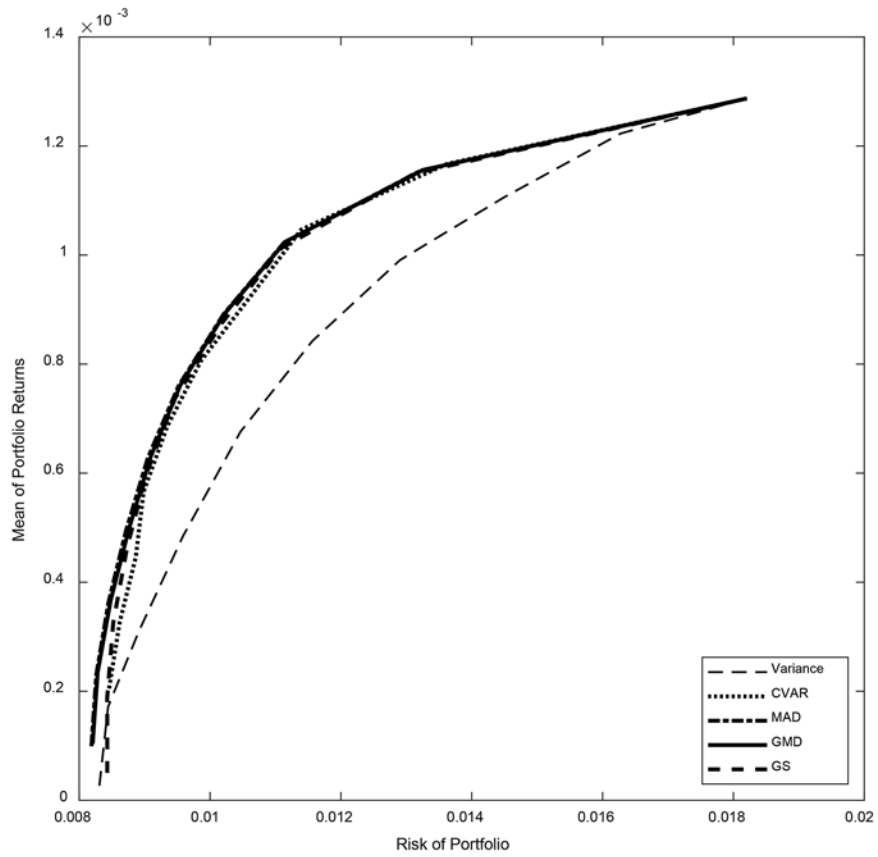

**Fig. 13.** Efficient frontiers EURO STOXX 50, simulation scenarios (VAR, CVAR, MAD, GMD, and $GS_2$)

As a starting point, we must observe that, as expected, MAD and GMD are very close in performance. This is clearly indicated in all the figures above where often the two measures are even undistinguishable. Our simulations with different stock selections produce the same pattern. There is robustness of this conclusion across different sample sizes, so it is not driven by numerical issues coming from consistency.

Now, moving to the performance of the measures, we can see that obviously the variance measure is distinct from the rest. This reflects its discussed suboptimality, both from the perspectives of variability and tail risk measurement. For other stock selections, the patterns are even more pronounced, with the winner-takes-all problem clearly visible in many cases. E.g., such is depicted very ostensibly both in Fig. 6 and Fig. 8. There the highest point on the efficient frontier for the variance case is achieved for a dominating stock which is not the same compared to its counterparts for the other measures where there is uniformity. For variance, appropriate diversification is often lacking – not only in the final portfolio on the efficient frontier but also in several portfolios preceding it. Now, we can proceed to more core observations.

The $GS_2$ measure is clearly different from the GMD and MAD ones and maintains some resemblance both to them and the CVAR measure. Its complex nature of combining variability and tail risk characteristics becomes evident. Furthermore, a careful look will reveal that the measure maintains not only a comparable degree of diversification to that offered by the other measures, but even more so at several points on the efficient frontier. We intentionally present cases where the effect is not distinctly clear-cut. For other stock selections, it is much stronger, explicitly highlighting the better diversification qualities of $GS_2$. Additionally, the choice of weights for the tails and their cut-off values in the construction of the measure exert influence. A careful calibration is necessary and can be postponed to future specialized portfolio analysis. Effectively, such will exactly position the measure between its two natures and their apparitions. It is exactly the $GS_2$ flexibility that makes it a plausible choice that might be considered by portfolio and risk managers.

We can further observe that we somehow artificially chose $v = 2.5$ for our exposition. Precise values for the coefficient at specific time points can be accurately estimated by following the logic[31] outlined in (*Shalit & Yitzhaki*, 1989). There exactly the one of 2.5 was found to be good for fitting the market the authors focus on. In the current context, not the specific market is of primary interest but the plausibility of the value for capturing risk aversion and thus its suitability for analysis. The latter is obviously valid, and we can work with risk estimates around the magnitude elaborated. Our simulations give that for larger values of υ, we get close results between $GMD(r_\alpha,v)$ and $GS_2(r_\alpha,v)$. This resembles the already discussed situation of narrow outcomes between $GMD(r_\alpha,v)$ and $GS_1(r_\alpha,v)$. For lower values of υ, but still within the risk aversion range, we get more diversified portfolio compositions for both risk measures. However, we get also larger disparities between them with the $GS_2(r_\alpha,v)$ measure producing more diversification. The effect is stronger the closer the coefficient to the boundary value of 1 is. The reason for this pattern is that by inspecting Fig. 5, we will observe more relative weight given to the scenarios falling in the upper right corner of the plot. This means more risk for tail events to which the investor naturally responds by more diversification.

## Conclusion

The setting from the main body proves useful for generating ideas in quantitative portfolio management. A specific set of risk measures was suggested, showing promising results. However, the analysis requires focused attention for a more in-depth and comprehensive treatment. This is especially valid for a proper calibration of the parameters characterizing the measures. Such would require the development of faster optimization routines resorting to linear programming. Further structural motivation for the proposed risk measures would also be of service.

[31] Within a representative agent framework, the risk preferences of the average market investor are indirectly revealed through the composition of a well-chosen market portfolio. This can be approximated by the value of a broad market index or, with plenty of caveats, by the equity index under focus. Using an appropriate metric, such as mean squared error, we can closely align the weights of this index with the weights of a portfolio from the efficient frontier of $GMD(r_\alpha,v)$ or $GS_2(r_\alpha,v)$ for fixed expected value of the returns.